\documentclass{amsart}
\usepackage[mathscr]{eucal}
\usepackage{amssymb, amsmath,array, amscd}
\usepackage[OT2,T1]{fontenc}

\newtheorem{thm}{Theorem}

\newtheorem{prop}{Proposition}

\newtheorem{lemma}{Lemma}[section]
\newtheorem{claim}{Claim}[section]

\newtheorem{definition}{Definition}

\newtheorem{rem}{Remark}[section]

\newtheorem{cor}{Corollary}[section]

\newtheorem{prob}{Problem}

\newtheorem{ex}{Example}

\def\limn{\underset{n\to\infty}{\ell im}}
\def\M{\mathcal{M}}

\def\V{\mathcal{V}}
\def\de{\mathcal{D}}

\def\pr{{\,\prime}}
\def\Ga{\Gamma}

\def\p{\mathbb{P}}
\def\wt{\widetilde}

\def\pa{\partial}
\def\sup{\supset}
\def\sub{\subset}

\def\I{\mathcal{I}}
\def\d{\delta}
\def\A{\mathcal{A}}
\def\emp{\emptyset}
\def\ov{\overline}
\def\om{\omega}
\def\Om{\Omega}
\def\fa{\mathcal{F}}
\def\a{\alpha}
\def\be{\beta}
\def\R{\mathbb{R}}
\def\ep{\epsilon}
\def\C{\mathbb{C}}

\def\{{\lbrace}
\def\}{\rbrace}

\def\la{\lambda}

\def\L{\mathcal{L}}
\def\g{\gamma}

\def\N{\mathbb{N}}
\def\Si{\Sigma}

\def\Z{\mathbb{Z}}
\def\O{\mathcal{O}}
\def\*{\star}

\def\r{\mathcal{R}}

\def\Q{\mathbb{Q}}

\def\var{\varphi}

\def\X{\mathcal{X}}

\def\wh{\widehat}

\begin{document}
\title[Commuting vector fields]{Commuting vector fields}


\author[D. Cerveau]{D. Cerveau}

\author[A. Lins Neto]{A. Lins Neto}

\thanks{The $1^{st}$ author was supported by CNRS and ANR-16-CE40-0008 project "Foliage".}
\thanks{The $2^{nd}$ author was partially supported by CNPq (Brazil) and University of Rennes 1.}



\begin{abstract}
In this paper we study the centralizer $C(X)$ of a germ of vector field at $0\in\C^n$.
A particular atention is given to the case of dimension two.
\end{abstract}

\maketitle

\tableofcontents

\section{Introduction}\label{ss:1}
Let $\O_n$ be the ring of germs at $0\in\C^n$ of holomorphic functions and $\X_n$ be the $\O_n$-modulus of germs at $0\in\C^n$ of holomorphic vector fields.
We say that $X,Y\in\X_n$ {\it commute} if $[X,Y]=0$, where $[.,.]$ denotes the Lie bracket.
We denote by $C(X)$ the set of germs of $\X_n$ commuting with $X$. Note that $C(X)$ is a $\C$-vector subspace of $\X_n$ and its complex dimension will be denoted as $d(X)$.
The purpose of this paper is to give examples and properties of the spaces $C(X)$ for a large class of vector fields $X$, in particular in dimension two.

For instance, if $\la\in\C$ then $\la X\in\C(X)$, so that $d(X)\ge1$.
When $d(X)=1$ then $C(X)=\C.\,X$ and we will say also that $C(X)$ is {\it trivial}.

There are examples in which
$d(X)=\infty$.
For instance, if $X$ has a non-constant holomorphic first integral, say $f$, and $g=\Phi\circ f$, where $\Phi\in\O_1$ then $Y=g.\,X\in C(X)$, because $g$ is a first integral of $X$:
\[
[X,Y]=[X,g.\,X]=X(g).\,X+g.\,[X,X]=0\,.
\]

In theorem \ref{t:fi} of \S\,\ref{ss:21} we will prove the converse when $n=2$: if $X\in\X_2$ and
$d(X)=\infty$ then $X$ has a non-constant holomorphic first integral.

Another important observation is that $C(X)$ is a Lie algebra: if $Y,Z\in C(X)$ then $[Y,Z]\in C(X)$, which is a consequence of Jacobi's identity
\[
[X,[Y,Z]]=[[X,Y],Z]+[Y,[X,Z]]=0\,.
\]
If $X$ has no non-constant meromorphic first integral then $C(X)$ is a finite dimensional Lie algebra (see proposition \ref{p:8} of \S\,\ref{ss:22}). 
\vskip.1in

Given $Y_1,...,Y_r\in\X_n$, we say that they are {\it generically linearly independent} (briefly g.l.i) if
$Y_1\wedge...\wedge Y_r\not\equiv0$. In other words, the analytic subset of $(\C^n,0)$ given by
\[
\{z\in(\C^n,0)\,|\,Y_1(z)\wedge...\wedge Y_r(z)=0\}
\]
is proper.
Given $X\in\X_n$ we define:
\[
r(X)=max\{r\in\N\,|\,\text{there are g.l.i}\,Y_1,...,Y_r\in C(X)\}\,.
\]

Note that, in general
\[
1\le r(X)\le\,min(n,d(X))\,.
\]
Moreover, if $d(X)>r(X)$ then $X$ has a non-constant meromorphic first integral (see proposition \ref{p:2} in \S\,\ref{ss:21}).

\vskip.1in
Let us see some simple examples.
\begin{ex}
{\rm Let $X=\frac{\pa}{\pa z_1}\in\X_n$, $n\ge2$. Then $Y\in C(X)$ if, and only if, $Y=f(z_2,...,z_n).\,v$, where $f\in\O_{n-1}$ and $v$ is a constant vector field.
We can say that
\[
C(X)\simeq\,\C\{z_2,...,z_n\}.\,\C^n\,,
\]
so that $r(X)=n$ and $d(X)=\infty$.}
\end{ex}

\begin{ex}
{\rm The radial vector field in $\C^n$ is given by $R=\sum_{j=1}^nz_j\frac{\pa}{\pa z_j}$.
It is easy to check that $Y\in\X_n$ commutes with $R$ if, and only if, $Y$ is a linear vector field
\[
Y=\sum_{1\le i,j\le n}a_{ij}\,z_i\,\frac{\pa}{\pa z_j}\,,
\]
where $(a_{ij})$ is a $n\times n$ matrix with constant entries.
In particular, we have $r(X)=n$ and $d(X)=n^2$.}
\end{ex}

\begin{ex}
{\rm Let $X$ be the diagonal vector field $X=\sum_{j=1}^n\la_j\,z_j\,\frac{\pa}{\pa z_j}$, where we will assume that $0\ne\la_i\ne\la_j$ $\forall$ $i\ne j$.
It is easy to see that if $Y\in C(X)$ is a linear vector field then $Y$ is also a diagonal vector field $Y=\sum_{j=1}^n\mu_j\,z_j\,\frac{\pa}{\pa z_j}$.
If the eigenvalues $\la_1,...,\la_n$ satisfy the non-resonant conditions:
\[
\la_i\ne \sum_{1\le j\le n}k_j\,\la_j\,,\,\forall\,1\le i\le n\,,\,\forall\,k_1,...,k_n\in\Z_{\ge\,0}\,\text{with}\,\sum_jk_j\ge2
\]
then $C(X)=$ the set of diagonal vector fields and $d(X)=r(X)=n$.

However, if the eigenvalues have a resonance, say $\la_i=\sum_{j=1}^nk_j\,\la_j$, then the non-linear vector field $Y=\Pi_{j=1}^nz_j^{k_j}.\,\frac{\pa}{\pa z_i}$ commutes with $X$.
If we denote the set of diagonal vector fields by $\de_n$ then
\[
C(X)=\de_n\,\oplus\,\left<\Pi_{j=1}^nz_j^{k_j}.\,\frac{\pa}{\pa z_i}\,|\,\la_i=\sum_jk_j\,\la_j\right>_\C\,,
\]
where the notation $\left<A\right>_\C$ denotes the $\C$-vector space generated by the set $A$.
In particular, if the eigenvalues have a resonance then $r(X)=n<d(X)\le\infty$.}
\end{ex}

More examples in the case $n=2$ will be seen in \S\,\ref{ss:22}.
\vskip.1in
In \S\,\ref{ss:3} the case of dimension two will be studied in detail.
We will see that in this case, if $r(X)=2$ then the foliation $\fa_X$, induced by $X$, is Liouvillian integrable: it can be defined by a closed meromorphic 1-form (see \S\,\ref{ss:31}).

In \S\,\ref{ss:32} we study $C(X)$ when $X\in\X_2$ is a generalized curve (see \cite{cls}) and has just one {\it separatrix}.
In this case, we will see that $r(X)=1$ and that, either $d(X)=1$, or $d(X)=\infty$ (if $X$ has a non-constant holomorphic first integral).
Recall that a separatrix of $X$ is a germ of curve through the origin, say $\g\colon(\C,0)\to(\C^n,0)$, which is regular outside $0$ and is $X$-invariant. This means that $\g(0)=0$ and $0\ne\g^\pr(t)\in\C.\,X(\g(t))$ if $t\ne0$.
In the two dimensional case we say that $X$ is {\it non-dicritical} if it has a finite number of irreducible separatrices. Otherwise, we say that it is {\it dicritical}.

In \S\,\ref{ss:24} we study the so-called generalized curves.
We will see that if $X$ is a non-dicritical generalized curve, with an isolated singularity at the origin, $DX(0)$ is nilpotent and has in its reduction of singularities a singularity with non-rational characteristic number then $C(X)$ is trivial.

In \S\,\ref{ss:34} we study homogeneous and quasi-homogeneous vector fields.
We will see that if $X\in\X_2$ is quasi-homogeneous and $r(X)=2$ then $X$ has a non-constant meromorphic first integral.
Another result, in this case, is that if $X$ is non-dicritical has an isolated singularity at $0\in\C^2$ and $DX(0)$ is nilpotent then $r(X)=1$ and, either $d(X)=1$, or $X$ has a non-constant holomorphic first integral and $d(X)=\infty$ (theorem \ref{t:3}).

In \S\,\ref{ss:dic} we study the dicritical case, when the vector field $X$ has infinitely many separatrices through $0\in\C^2$. We will see that $1\le d(X)\le4$ and we will describe completely the cases $d(X)=3$ and $d(X)=4$:
if $d(X)\in\{3,4\}$ then $X$ is linearizable and, modulo a multiplicative constant, $X$ is conjugated to the radial vector field (if $d(X)=4$) or to a linear vector field of the form $z_1\frac{\pa}{\pa z_1}+n\,z_2\frac{\pa}{\pa z_2}$, where $n\in\N_{\ge2}$ (if $d(X)=3$).
When $d(X)\le 2$ we will see that $d(X)=r(X)$ in the dicritical case.
\vskip.1in
We finish this section by fixing more notations that be will used along the paper.
\begin{itemize}
\item[1.] 
$\wh{\O}_n=$ the ring of formal power le  at $0\in\C^n$.

$\M_n=\M(\C^n,0)=$ the field of germs at $0\in\C^n$ of meromorphic functions.
\item[2.] $Diff(\C^n,0)=$ the group of germs of holomorphic diffeomorphisms of $(\C^n,0)$ fixing $0\in\C^n$.
\item[3.]
Given a subring $\r\sub\O_n$ and vector fields $Y_1,...,Y_s\in \X_n$
we denote
\[
\left<Y_1,...,Y_s\right>_\r:=\{a_1.\,Y_1+...+a_s.\,Y_s\in\X_n\,|\,a_1,...,a_s\in\r\}
\]
the sub $\r$-modulus generated by the $Y_{j`s}$.
\item[4.] Given a germ $Y$ at $0\in\C^n$ of holomorphic vector field, function or form, we will denote by $\mu(Y,0)$ its algebraic multiplicity at $0$ (the order of the first non-zero jet of $Y$).\item[5.] The ring of holomorphic first integrals of a vector field $X\in\X_n$ will be denoted by
$\I(X)$: $\I(X)=\{f\in\O_n\,|\,X(f)=0\}$.
\end{itemize}

If $z=(z_1,...,z_n)$ is a local coordinate system, and $X=\sum_{i=1}^nX_i(z)\frac{\pa}{\pa z_i}$, then $f\in\I(X)$ iff
\[
X(f):=\sum_{i=1}^nX_i(z)\,\frac{\pa f}{\pa z_i}=0\,.
\]

The field of meromorphic first integrals of $X\in\X_n$ will be denoted by $\M(X)$:
\[
\M(X)=\{f\in\M_n\,|\,X(f)=0\}
\]
Note that $\C\sub\I(X)\sub\M(X)$.
A function $f\in\M(X)$ is said to be a {\it pure} meromorphic first integral if $f,1/f\notin\I(X)$.
\vskip.1in
Observe that $\I(X).\,X\sub\,C(X)$: if $f\in\I(X)$ then
\[
[X,f.\,X]=X(f).\,X+f.\,[X,X]=0\,.
\]
\vskip.1in
For instance, if $X(0)\ne0$, i.e. $X$ is non-singular, then in some local coordinate system $z=(z_1,...,z_n)$ we have $X=\frac{\pa}{\pa z_1}$.
In this case $\I(X)=\C\{z_2,...,z_n\}$, the ring of convergent power series on the variables $z_2,...,z_n$, and
\[
C(X)=\left<\frac{\pa}{\pa z_1},...,\frac{\pa}{\pa z_n}\right>_{\I(X)}\,.
\]

\begin{rem}
{\rm We have seen that $1\le r(X)\le min\{n,d(X)\}$.
About this inequality, we would like to remark the following:
\begin{itemize}
\item[a.] When $r(X)<d(X)$ then $X$ has a non-constant meromorphic first integral.
\item[b.] When $n=2$, $r(X)<d(X)$ and $0$ is an isolated singularity of $X$ then $\I(X)\supsetneq\C$, that is $X$ has a non-constant holomorphic first integral.
\item[c.] When $r(X)=n$, then $X$ is Liouvillian integrable: there exists a closed meromorphic $(n-1)$-form $\om$ such that $i_X\om=0$.
\end{itemize}
Some of the above remarks will be proved in \S\,\ref{ss:2}.}
\end{rem}

\section{First properties and examples}\label{ss:2}

\subsection{General properties.}\label{ss:21}
In this section we state some elementary properties that will be used along the paper.
Some of the remarks stated in \S\,\ref{ss:1} will be proved here.

\begin{prop}\label{p:1}
Let $X\in\X_n$, $X\not\equiv 0$, with $r(X)=r\ge1$. Given $z\in(\C^n,0)$ consider the vector subspace $\de(X)(z)$ of $\C^n$ defined as
\[
\de(X)(z):=\{Y(z)\in\C^n\,|\,Y\in C(X)\}\,.
\]
Then:
\begin{itemize}
\item[(a).] The set $V:=\{z\in(\C^n,0)\,|\,dim_\C(\de(X)(z))=r\}$ is the complement of a proper analytic subset on $(\C^n,0)$. In particular, it is open and dense in $(\C^n,0)$.
\item[(b).] $\de(X)$ defines an integrable distribution of dimension $r$ on $V$.
\end{itemize}
\end{prop}

{\it Proof.} For simplicity, from now on we will work with representatives of the germs (of functions or vector fields) without specify their domain: we will write $z\in(\C^n,0)$ to denote that $z$ belongs to some domain where some representative of the germs are defined.

Consider the set
\[
S:=(\C^n,0)\setminus V=\{z\in (\C^n,0)\,|\,dim_\C(\de(X)(z))<r\}\,.
\]
Note that $z\in S$ if, and only if, for any $Y_1,...,Y_r\in C(X)$ then $Y_1(z)\wedge...\wedge Y_r(z)=0$.
In particular, we have
\[
S=\bigcap_{\underset{1\le j\le r}{Y_j\in C(X)}}(Y_1\wedge...\wedge Y_r=0)\,\,.
\]
The set $S$ is analytic because it is an intersection of analytic subsets of $(\C^n,0)$.
On the other hand, by the definition of $r=r(X)$ there exist $Y_1,...,Y_r\in C(X)$ such that $Y_1\wedge...\wedge Y_r\not\equiv0$.
Hence $S$ is a proper analytic subset of $(\C^n,0)$. This proves (a).
\vskip.1in
Assertion (b) follows from the fact that $C(X)$ is a Lie algebra: $Y,Z\in C(X)$ $\implies$
$[Y,Z]\in C(X)$.
\qed
\vskip.1in
\begin{prop}\label{p:2}
Let $X\in\X_n$, $X\not\equiv0$, be such that $r=r(X)< d(X)$. Then $\M(X)\ne\C$.
\end{prop}

{\it Proof.} Let $Y_1,...,Y_r\in C(X)$ be generically independent vector fields of $C(X)$.
Since $d(X)>r$ there exists $Y\in C(X)$ such that $Y\notin\left<Y_1,...,Y_r\right>_\C$.
By the definition of $r=r(X)$ we have $Y\wedge Y_1\wedge...\wedge Y_r\equiv0$.
Let
\[
V:=\{z\in(\C^n,0)\,|\,Y_1(z)\wedge...\wedge Y_r(z)\ne0\}\,.
\]
Note that $V$ is the complement of an analytic subset of $(\C^n,0)$.
Moreover, if $z\in V$ then $Y(z)\in\left<Y_1(z),...,Y_r(z)\right>_\C$, so that we can write
\begin{equation}\label{eq:Y}
Y(z)=f_1(z).\,Y_1(z)+...+f_r(z).\,Y_r(z)\,,
\end{equation}
where $f_1,...,f_r\in\O(V)$. Since
\[
Y_1\wedge...\wedge Y_{j-1}\wedge Y\wedge Y_{j+1}\wedge...\wedge Y_r=f_j.\,Y_1\wedge...\wedge Y_r\,,\,1\le j\le r\,,
\]
the functions $f_1,...,f_r$ extend to meromorphic functions in a neighborhood of $0$: $f_j\in\M_n$, $1\le j\le r$.
Since $Y\notin\left<Y_1,...,Y_r\right>_\C$ at least one of the functions is non-constant.
Finally, relation (\ref{eq:Y}) implies
\[
\sum_{j=1}^rX(f_j).\,Y_j=[X,Y]=0\,\implies\,X(f_j)=0\,,\,1\le j\le r\,\,.\qed
\]

\begin{cor}\label{c:21}
Let $X\in\X_n$ having an isolated singularity at $0\in\C^n$. Then:
\begin{itemize}
\item[1.] If $r(X)=1<d(X)$ then $\I(X)\ne\C$.
\item[2.] If $r(X)<d(X)$ then $\M(X)\ne\C$. In particular, if $n=2$ and $X$ is non-dicritical (has finitely many separatrices through the origin) then $\I(X)\ne\C$.
\end{itemize}
\end{cor}
{\it Proof.}
If $r=r(X)=1$ and $d(X)>1$, then there exists $Y\in C(X)\setminus\C.\,X$. Note that $Y\wedge X=0$, because $r=1$. Since $X$ has an isolated singularity at $0$ we have $Y=f.\,X$, where $f\in\O_n\setminus\C$. Finally, $X(f)=0$ and so $f\in\I(X)\setminus\C$.
\vskip.1in
If $r(X)<d(X)$ then, as we have seen in proposition \ref{p:2}, $X$ has a non-constant meromorphic first integral, say $f\in\M(X)$. If $n=2$ and $f$ is purely meromorphic then $X$ is dicritical.
\qed
\vskip.1in
In the next section we will see examples of $X\in\X_2$ such that $r(X)=2<d(X)$ and $\I(X)=\C$. By corollary \ref{c:21} these examples are dicritical.
\vskip.1in 
When $\M(X)\ne\C$ then $d(X)$, although finite, can be arbitrarily big, as shows the following example in any dimension $n\ge2$:
\begin{ex}\label{ex:5}
{\rm Let $R=\sum_{j=1}^nz_j\frac{\pa}{\pa z_j}$ be the radial vector field and $X=(z_1...z_n)^k.\,R$, where $k\in\N$. It can be verified that
$C(X)=\L\,\oplus\,\A$, where $\L=$
\[
=\left\{L\,|\,L=\sum_{j=1}^n\la_j\,z_j\frac{\pa}{\pa z_j}\,\,\text{and}\,\,\sum_{j=1}^n\la_j=0\right\}=\{L\in\X_n\,|\,L\,\text{is linear and}\,\,L(z_1...z_n)=0\}
\]
and
\[
\A=\{g.\,R\,|\,g\,\,\text{is a homogeneous polynomial of degree $k.n$}\}\,.
\]}
\end{ex}

Example \ref{ex:5} motivates the following question:
\begin{prob}\label{pro:1}
Let $X\in\X_n$, where $X\not\equiv0$. Does $d(X)=\infty$ implies that $\I(X)\ne\C$?
\end{prob}
In dimension two problem \ref{pro:1} has a positive answer:
\begin{thm}\label{t:fi}
Let $X\in X_2$, $X\not\equiv0$, be a vector field with $d(X)=\infty$. Then $\I(X)\ne\C$.
\end{thm}

{\it Proof.}
We will assume first that the vector space $D(X)=\{Y\in C(X)\,|\,Y\wedge X=0\}$ is infinite dimensional.

By corollary \ref{c:21} we can assume that the origin is not an isolated singularity of $X$.
In particular, we can write $X=g.\,T$, where $g\in\O_2$, $g(0)=0$ and $T\in\X_2$ has an isolated singularity at $0\in\C^2$.

If $Y\in D(X)\setminus\{0\}$ then $Y\wedge T=0$ and since $0$ is an isolated singularity of $T$ we can write $Y=f_Y.\,T$, where $f_Y\in\O_2\setminus\{0\}$.
Note that $f_Y/g\in\M(X)$ because $Y\in C(X)$ and $Y=\frac{f_Y}{g}\,X$. In particular, we can assume that $f_Y(0)=0$ for all $Y\in D(X)$, for otherwise $g/f_Y\in\I(X)\setminus\C$ and we are done.
From now on, we will assume that $\M(X)\ne\C$ and, by contradiction, that $\I(X)=\C$.

Denote $\V(X):=\left<f_Y\,|\,Y\in D(X)\right>_\C$, the $\C$-vector subspace of $\O_2$ generated by the $f_Y$, $Y\in D(X)$. Since $dim(D(X))=\infty$ we must have $dim(\V(X))=\infty$.
As a consequence, the set $\mu(X):=\{\mu(f_Y,0)\,|\,Y\in D(X)\}\sub\N$ is unbounded.

In fact, let $j^N\colon\O_2\to\O_2/m_2^{N+1}$ be the $N^{th}$-jet map and $J^N=j^N|_{\V(X)}$ be its restriction to $\V(X)$. Note that $ker(J^N)=\left<f_Y\in\V(X)\,|\,\mu(f_Y,0)>N\right>_\C$. 
Since $dim(J^N(\V(X)))<\infty$ we get $dim(ker(J^N))=\infty$ for all $N\in\N$.
Hence, $\mu(X)$ is unbounded.

Let $\{Y_n\}_n$ be a sequence of vector fields in $D(X)$ such that the sequence
$\mu_n:=\mu(f_{Y_n},0)$ is increasing.

Recall that, by Seidenberg`s resolution theorem \cite{se}, there exists a blowing-up process
$\Pi\colon(M,E)\to(\C^2,0)$ such that the strict transform $\Pi^*(\fa_T)$ of the foliation $\fa_T$ is generically transverse to some irreducible component $D$ of the exceptional divisor $E$. In particular, at some generic point $p\in D$ there are local coordinates $(x,t)\colon(M,p)\to(\C^2,0)$ such that:
\begin{itemize}
\item[(i).] $x(p)=t(p)=0$ and $D_p=(x=0)$.
\item[(ii).] $\Pi^*(\fa_X)$ is defined at $p$ by $dt=0$.
\item[(iii).] $g\circ\Pi(x,t)=x^k.\,\wt{g}(x,t)$, where $x\not|\,\wt{g}$. We can choose $p$ in such a way that $\wt{g}$ is an unity.
\item[(iv).] $f_{Y_n}\circ\Pi(x,t)=x^{k_n}.\,\wt{f}_n(x,t)$, where $x\not|\,\wt{f}_n$.
\end{itemize}
Set $h_n:=\frac{f_{Y_n}}{g}\circ\Pi=x^{k_n-k}.\,\wt{h}_n$, where $\wt{h}_n=\wt{f}_n/\wt{g}$. Since $f_{Y_n}/g\in\M(T)$ we get $h_n\in\M(dt)$ and this implies
\[
dh_n\wedge dt=0\,\implies\,h_n=h_n(t)\,\implies\,k_n=k\,\,\forall\,n\in\N\,.
\]
However, it is known that if $f_1,f_2\in\O_2$ and $\mu(f_1,0)<\mu(f_2,0)$ then $\mu(f_1\circ\Pi,p)<\mu(f_2\circ\Pi,p)$ and this implies that $\limn k_n=\infty$, a contradiction.
\vskip.1in
Let us consider the general case: $d(X)=\infty$. We will assume also that $dim(D(X))<\infty$, so that $dim(C(X)/D(X)=\infty$.
Fix a local coordinate system $(x,y)$ around $0\in\C^2$ and set
$\mu=\frac{\pa}{\pa x}\wedge \frac{\pa}{\pa y}$.

Given $Y\in C(X)\setminus D(X)$ we can write $Y\wedge X=f_Y.\,\mu$, where $f_Y\in \O_2$.
Note that
\[
L_X\,Y\wedge X=[X,Y]\wedge X+Y\wedge [X,X]=0\,\implies
\]
\[
L_X(f_Y\,\mu)=(X(f_Y)\,-\nabla\,X.\,f_Y)\mu=0\,\implies
\]
\begin{equation}\label{eq:na}
\,\frac{X(f_Y)}{f_Y}=\nabla\,X\,.
\end{equation}
where $\nabla\,X=\frac{\pa X(x)}{\pa x}+\frac{\pa X(y)}{\pa y}$.
Let $\{Y_n\}_n$ be a sequence in $C(X)\setminus D(X)$ such that the sequence $\mu_n:=\mu(f_{Y_n},0)$ is increasing and set $f_n:=f_{Y_n}$.
From (\ref{eq:na}) we get
\[
\frac{X(f_n/f_1)}{f_n/f_1}=\frac{X(f_n)}{f_n}-\frac{X(f_1)}{f_1}=0\,\implies\,\frac{f_n}{f_1}\in\M(X)\,,\,\forall\,n\ge2\,.
\]
Therefore, we can apply the same argument of the first case to the sequence of meromorphic first integrals $f_n/f_1$.
\qed
\vskip.1in
We finish this section with the following result:
\begin{prop}\label{p:3}
Let $X\in\X_n$, $X\not\equiv0$. If $r(X)=n$ then there exists a closed meromorphic $(n-1)$-form anihilating $X$.
\end{prop}

{\it Proof.}
Fix a local coordinate system $(z_1,...,z_n)$ around $0\in\C^n$ and let $\nu=dz_1\wedge...\wedge dz_n$ and $X=\sum_{j=1}^nX_j\frac{\pa}{\pa z_j}$.
Set
\begin{equation}\label{eq:om0}
\om=i_X\nu=\sum_{j=1}^n(-1)^{j+1}X_j\,dz_1\wedge...\wedge\wh{dz_j}\wedge...\wedge dz_n\,,
\end{equation}
where $\wh{dz_j}$ means the omission of $dz_j$ in the product.
Since $r(X)=n$ there are vector fields $Y_1,...,Y_n\in C(X)$ such that $Y_1\wedge...\wedge Y_n\not\equiv0$.
Set 
\begin{equation}\label{eq:mu}
Y_1\wedge...\wedge Y_n:=g.\,\mu\,,
\end{equation}
where $g\in\O_n$ and $\mu=\frac{\pa}{\pa z_1}\wedge...\wedge\frac{\pa}{\pa z_n}$.
From (\ref{eq:om0}) we get
\[
d\,\om=\sum_j(-1)^{j+1}dX_j\wedge dz_1\wedge...\wedge\wh{dz_j}\wedge...\wedge dz_n=\nabla X.\,\nu\,,
\]
where $\nabla X=\sum_j\frac{\pa X_j}{\pa z_j}$. Since $Y_j\in C(X)$, $1\le j\le n$, we get from (\ref{eq:mu})
\[
0=L_X\left(Y_1\wedge...\wedge Y_n\right)=L_X(g.\,\mu)=X(g).\,\mu+g.\,L_X\mu=X(g).\,\mu-\nabla X.\,g.\,\mu\,\implies
\]
$X(g)=\nabla X.\,g$.
In particular
\[
g.\,d\,\om=g.\,\nabla X.\,\nu=X(g).\,\nu=dg\wedge i_X\nu=dg\wedge\om\,\,\implies\,\,d\left(\frac{\om}{g}\right)=0\,\,.
\]
Therefore, the meromorphic form $\frac{\om}{g}$ is closed and satisfies $i_X\frac{\om}{g}=0$.
\qed

\subsection{More general remarks.}\label{ss:22a}

We begin by the case in which $C(X)$ is a finite dimensional Lie algebra. 

\begin{prop}\label{p:8}
Let $X\in\X_n$ be a germ of holomorphic vector field, where $n\ge2$. Suppose that $\M(X)=\C$, that is $X$ has no non constant meromorphic first integral. Then $C(X)$ is a finite dimensional Lie algebra of dimension $r(X)$. In particular, $d(X)=r(X)$.
\end{prop}

{\it Proof.}
Let $Y_1,...,Y_r\,\in C(X)$, $r=r(X)$, be such that $Y_1\wedge...\wedge Y_r\not\equiv0$.
If $Y\in C(X)$ then there exist $f_1,...,f_r\in \M_n$ such that
\[
Y=\sum_{j=1}^rf_j.\,Y_j\,.
\]
In particular, we have:
\[
0=[X,Y]=\sum_{j=1}^rX(f_j)\,Y_j\,\implies\,X(f_j)=0\,,\,1\le j\le r\,\,.
\]
Since $\M(X)=\C$ we must have $f_j\in\C$, $1\le j\le r$.
\qed

\vskip.1in
A particular case of proposition \ref{p:8} is the following:
\begin{cor}\label{c:41}
Let $X\in\X_n$, where $n\ge2$. Suppose that $\M(X)=\C$ and that $C(X)\sub m_n^2.\,\X_n$. Then $C(X)$ is nilpotent.
\end{cor}

Corollary \ref{c:41} is a direct consequence of the following:
\begin{prop}\label{p:9}
Let $\L\sub\X_n$ be a finite dimensional Lie algebra of germs of vector fields.
Suppose that $\L\sub m_n^2.\,\X_n$. Then $\L$ is nilpotent. 
\end{prop}

{\it Proof.} Proposition \ref{p:9} is a direct consequence of Engel-Lie theorem \cite{hu}.
In fact, given $Y\in\L$ consider the operator
\[
ad_Y\colon\L\to\L\,,\,ad_Y(Z)=[Z,Y]\,\,.
\]
Then the eigenvalues of $ad_Y$ are all zero, because
\[
\mu([Z,Y],0)\ge\mu(Z,0)+\mu(Y,0)-1>\mu(Z,0)\,.
\]
Therefore $\L$ is nilpotent by the Engel-Lie theorem \cite{hu}.
\qed
\vskip.1in
\begin{rem}
{\rm In fact, proposition \ref{p:9} can be generalized to finite dimensional Lie algebras $\L\sub m_n.\,\X_n$ such that the linear part $j^1\,\L$ is nilpotent.}
\end{rem}
As a consequence we have:
\begin{cor}
Let $\L\sub m_n.\,\X_n$ be a finite Lie algebra and $J^1\colon\L\to gl(n,\C)$ be the morphism $J^1(X)=DX(0)$. Then $ker(J^1)=\L\cap m_n^2.\,\X_n$ is nilpotent.
\end{cor}

\begin{rem}
{\rm In example \ref{ex:5} $C(X)=\L\,\oplus\,\A$ is solvable but not nilpotent: in fact, as the reader can check $[\L\,\oplus\,\A,\L\,\oplus\,\A]=[\L,A]\sub\A$, and $[\A,\A]=0$, but $[\L,[\L,\A]]=[\L,\A]$.}
\end{rem}

Another important result in this direction is the following (cf. \cite{he} and \cite{gs}):

\begin{thm}
{\rm (R. Hermann, V. Guillemin, S. Sternberg)}
Let $\L\sub m_n.\,\X_n$ be a semi-simple Lie algebra. Then $\L$ is linearizable, that is there exists $\phi\in Diff(\C^n,0)$ such that $\phi_*(\L)$ is a Lie algebra of linear vector fields.
\end{thm}

As a consequence, we have the following:

\begin{cor}
Let $X\in m_n.\,\X_n$ and assume that $\M(X)=\C$. Then, up to conjugacy the Levi-Malcev decomposition {\rm(see \cite{jn})} is of the form $C(X)=Rad\,\oplus\,\L$, where $\L$ is a semi-simple Lie algebra of linear vector fields and $Rad$ is the solvable radical of $C(X)$. 
\end{cor}

Some easy consequences in small dimension when $\M(X)=\C$:

\begin{itemize}
\item[1.] If $r(X)=1$ then $C(X)=\C.\,X$.
\item[2.] If $r(X)=2$ then $C(X)$ is abelian.
\item[3.] If $r(X)=3$ then, either $C(X)$ is abelian, or $C(X)\simeq\C\,\bigoplus\,A$, where $A$ is the affine Lie algebra, or $C(X)$ is isomorphic to the Heisenberg algebra:
\[
C(X)=\left<X,Y,Z\right>_\C\,\,\text{where}\,\,[X,Y]=[X,Z]=0\,\,\text{and}\,\,[Y,Z]=X\,\,.
\]
\end{itemize}

\begin{ex}
{\rm Let $X$ be the "Poincar\'e-Dulac" vector field defined as
\[
X=(m\,x+y^m)\frac{\pa}{\pa x}+y\frac{\pa}{\pa y}+(n\,z+y^n)\frac{\pa}{\pa z}\,,\,\text{where}\,m,\,n\in\N\,.
\]
Then $\M(X)=\C$ and $C(X)=\left<y^m\frac{\pa}{\pa x},y^n\frac{\pa}{\pa z},X\right>_\C$ is abelian of dimension three.}
\end{ex}

\begin{prob}
Is there $X\in\X_3$ such that $\M(X)=\C$, $r(X)=3$ the origin is an isolated singularity of $X$ and $C(X)$ is isomorphic to, either $\C\oplus\,A$ or to the Heisenberg algebra?
\end{prob}

In his work on commuting vector fields \cite{ln} the second author shows the following result:
\vskip.1in
{\bf Theorem \cite{ln}.}
{\it Let $X\in\wh{\X}_2$ with $X=\sum_{j\ge k}X_j$ its Taylor series, where $X_j$ is homogeneous of degree $j$ and $X_k\ne0$. Assume that 
\begin{itemize}
\item[(a).] $k\ge2$ and $X_k$ has an isolated singularity at $0\in\C^2$.
\item[(b).] $X_k$ has no meromorphic first integral.
\end{itemize}
Then $C(X)=\C.\,X$.}
\vskip.1in
A generalization of this result is the following:

\begin{prop}\label{p:com}
Let $X\in\X_n$ with $X=\sum_{j\ge k}X_j$ its Taylor series,  where $X_j$ is homogeneous of degree $j$ and $X_k\ne0$. Assume that 
\begin{itemize}
\item[(a).] $k\ge2$ and $X_k$ has an isolated singularity at $0\in\C^2$.
\item[(b).] The foliation $\fa_{X_k}$, induced by $X_k$ on $\p^{n-1}$, has no algebraic invariant set of dimension $\ell$, where $0<\ell<n-1$.
\end{itemize}
Then $C(X)=\C.\,X$.
\end{prop}

{\it Proof.}
The proof is based in the following:
\begin{lemma}\label{l:com}
Let $X_k$ be as in (a) and (b) of proposition \ref{p:com}. Let $Y_\ell\in\X_n$ be homogeneous of degree $\ell$ be such that $[X_k,Y_\ell]=0$. Then:
\begin{itemize}
\item[(i).] If $\ell\ne k$ then $Y_\ell=0$.
\item[(ii).] If $\ell=k$ then $Y_k=\la.\,X_k$, where $\la\in\C$.
\end{itemize}
\end{lemma}

Let us prove the proposition using lemma \ref{l:com}. Let $Y\in\X_n$ be such that $[X,Y]=0$ with Taylor series $Y=\sum_{j\ge \ell}Y_j$, where $Y_j$ is homogeneous of degree $j$ and $Y_\ell\ne0$.

Since $[X_r,Y_s]$ is homogeneous of degree $r+s-1$, we can write the Taylor series of $[X,Y]$ as
$[X,Y]=\sum_{j\ge k+\ell-1}Z_j$, where $Z_{k+\ell-1}=[X_k,Y_\ell]$.
Therefore $[X_k,Y_\ell]=0$, and by lemma \ref{l:com} we get $\ell=k$ and $Y_k=\la.\,X_k$, where $\la\in\C^*$.
Now, if we set $W=Y-\la.\,X$ then $[X,W]=0$ and
\[
W=\sum_{j\ge k+1}(Y_j-\la.\,X_j)\,\,\implies\,\,Y_j=\la.\,X_j\,\,,\,\,\forall\,j\,,
\]
by lemma \ref{l:com}. Hence, $Y=\la.\,X$ and $Y\in\C.\,X$.
\qed
\vskip.1in
{\it Proof of lemma \ref{l:com}.}
First of all, hypothesis (b) of the proposition implies that if $A$ is an irreducible analytic subset of $(\C^n,0)$ of dimension $dim(A)\ge1$ and $X_k$-invariant then $A$ is a straigth line through the origin, corresponding to a singularity in $\p^{n-1}$ of $\fa_{X_k}$.
In particular, we have $\M(X_k)=\C$.

Let $Y_\ell\ne0$ be homogeneous of degree $\ell$ such that $[X_k,Y_\ell]=0$.
It is sufficient to prove that $X_k\wedge Y_\ell\equiv0$. In fact, since $X_k$ has an isolated singularity at the origin this implies that $Y_\ell=f.\,X_k$, where $X_k(f)=0$. Since $\M(X_k)=\C$ we get $\ell=k$ and $Y_\ell\in\C.\,X_k$.

Suppose by contradiction that $X_k\wedge Y_\ell\ne0$.
Let $R=\sum_{j=1}^nz_j\frac{\pa}{\pa z_j}$ be the radial vector field.
Note that $X_k\wedge R\ne0$, because otherwise $0$ would not be an isolated singularity of $X_k$.
We have two possibilities:
\vskip.1in
$1^{st}$. $X_k\wedge Y_\ell\wedge R\equiv0$.
In this case, we can write
\begin{equation}\label{e:4}
Y_\ell=f.\,X_k+g.\,R\,\,,
\end{equation}
where $f$ and $g$ are meromorphic.
Note that $g\ne0$. In fact, from (\ref{e:4}) we get
\[
0\not\equiv Y_\ell\wedge X_k=g.\,R\wedge X_k\,\,\implies\,\,g\ne0\,.
\]
Again by (\ref{e:4}) we have
\[
0=[X_k,Y_\ell]=X_k(f).\,X_k+X_k(g).\,R+g.\,[X_k,R]=
\]
\[
=(X_k(f)-(k-1)\,g).\,X_k+X_k(g).\,R\,\,\implies\,\,X_k(g)=0\,\,\text{and}\,\,X_k(f)=(k-1)\,g\,\,,
\]
where in the above relation we have used that $[R,X_k]=(k-1).\,X_k$.

Since $\M(X_k)=\C$ we get $g\in\C^*$ and so $X_k(f)=(k-1)\,g\in\C^*$. Finally $X_k(f)=0$, because $0$ is a singularity of $X_k$, a contradiction with $g\ne0$.
\vskip.1in
$2^{nd}$. $X_k\wedge Y_\ell\wedge R\not\equiv0$. Let $A=\{z\in\C^n\,|\,X_k(z)\wedge Y_\ell(z)\wedge R(z)=0\}$.
Observe first that $dim(A)\ge2$: it is well known that the set of zeroes of a totally decomposable $r$-vector has dimension $\ge r-1$, unless it is empty.
We assert that $A$ is $X_k$-invariant.

In fact, first of all we have
\begin{equation}\label{e:5}
L_{X_k}\,(X_k\wedge Y_\ell\wedge R)=X_k\wedge Y_\ell\wedge[X_k,R]=0\,.
\end{equation}
If we write
\[
X_k\wedge Y_\ell\wedge R=\sum_{1\le i<j<k\le n}A_{ijk}(z)\,\frac{\pa}{\pa z_i}\wedge \frac{\pa}{\pa z_j}\wedge\frac{\pa}{\pa z_k}
\]
then the ideal defining $A$ is $\I(A)=\left<A_{ijk}\,|\,1\le i<j<k\le n\right>$.
Finally, relation (\ref{e:5}) implies that $X_k(\I(A))\sub\I(A)$, as the reader can check. Hence, $A$ is $X_k$-invariant and $dim(A)\ge2$, a contradiction with hypothesis (b).

Therefore, $X_k\wedge Y_\ell=0$ which proves the lemma.
\qed
\vskip.1in
\subsection{Examples.}\label{ss:22}
The aim of this section is to introduce some simple examples.

\subsubsection{Linear vector fields on $\C^2$ and the saddle-node.}
A non-zero linear vector field $X$ on $\C^2$ is one that can be written as $X=A(x,y)\frac{\pa}{\pa x}+B(x,y)\frac{\pa}{\pa y}$, where $A$ and $B$ are linear.
According to Jordan's normal form, after a linear change of variables and multiplication by a constant, it can be written in one of the following forms:

a. $X=X_\la=x\frac{\pa}{\pa x}+\la.\,y\frac{\pa}{\pa y}$, where $\la\in\C$ (semi-simple case).

b. $X=x\frac{\pa}{\pa y}$ (nilpotent case).

c. $X=x\frac{\pa}{\pa x}+(x+y)\frac{\pa}{\pa y}$.
\vskip.1in
The saddle-node is more complicated. The germ $X\in\X_2$ has a saddle-node at $0\in\C^2$ when it has an isolated singularity at $0$ and the linear part $DX(0)$ has one eigenvalue zero and the other non-zero. In this case, it is known that, after a {\it formal} change of variables and a multiplication by a constant, then the saddle-node can written as:

d. $X=x^{p+1}\,\frac{\pa}{\pa x}+y\,(1+\la\,x^p)\,\frac{\pa}{\pa y}$\,,

where $p\in\N$ and $\la\in\C$.
The formal invariants are the multiplicity $\mu=p+1$ and
$\la\in\C$.
The saddle-node in the formal normal form (d) is Liouvillian integrable, in the sense that the associated foliation can be defined also by a meromorphic closed form. The dual form, $\om:=i_Xdx\wedge dy$, has an integrating factor $d\left(\frac{1}{y.\,x^{p+1}}.\,\om\right)=0$ :  
\[
\frac{\om}{y.\,x^{p+1}}=\frac{1}{y.\,x^{p+1}}\,(x^{p+1}\,dy-y(1+\la\,x^p)\,dx)
=\frac{dy}{y}-\frac{dx}{x^{p+1}}-\la\frac{dx}{x}\,\,.
\]
For this reason, $\la$ is called the residue.
\vskip.1in
In the next table we specify $C(X)$, $r(X)$ and $d(X)$ in the above cases.
\vskip.1in
\begin{tabular}{|l||c|c|c|r|}\hline
$\,$&$\,\,\,\,X$&$C(X)$&$r(X)$&$d(X)$\\\hline
1.&$X_1=x\frac{\pa}{\pa x}+y\frac{\pa}{\pa y}$&all linear vector fields&2&4\\\hline
2.&$X_\la\,,\, \la\notin\Q_{\le0}\cup\N\cup 1/\N$&$\C.\,\left<x\,\frac{\pa}{\pa x},y\,\frac{\pa}{\pa y}\right>$&2&2\\\hline
3.&$X_{-p/q}\,,\,p,q\in\N$&$\C\{x^p.\,y^q\}.\,\left<x\,\frac{\pa}{\pa x},y\,\frac{\pa}{\pa y}\right>$&2&$\infty$\\\hline
4.&$X_0=x\,\frac{\pa}{\pa x}$&$\C\{y\}.\,\left<x\,\frac{\pa}{\pa x},y\,\frac{\pa}{\pa y}\right>$&2&$\infty$\\\hline
5.&$X_n\,,\,n\ge2$&$\C.\left<x\frac{\pa}{\pa x},y\frac{\pa}{\pa y},x^n\frac{\pa}{\pa y}\right>$&2&3\\\hline
6.&$x\,\frac{\pa}{\pa y}$&$\C\{x\}.\left<x\frac{\pa}{\pa x}+y\frac{\pa}{\pa y},\frac{\pa}{\pa y}\right>$&2&$\infty$\\\hline
7.&$x\,\frac{\pa}{\pa x}+(y+x)\frac{\pa}{\pa y}$&$\C.\left<x\frac{\pa}{\pa x}+y\frac{\pa}{\pa y},x\frac{\pa}{\pa y}\right>$&2&2\\\hline
8.&$x^{p+1}\frac{\pa}{\pa x}+y(1+\la\,x^p)\frac{\pa}{\pa y}$&$\C.\left<x^p\left(x\frac{\pa}{\pa x}+\la\,y\frac{\pa}{\pa y}\right),y\,\frac{\pa}{\pa y}\right>$&2&2\\\hline
\end{tabular}
\vskip.1in
In the above table, examples 1, 3, 4, 5 and 6 have $d(X)>r(X)$. In examples 3, 4 and 6 we have $\I(X)\ne\C$.
Examples 1 and 5 are dicritical and have purely meromorphic first integrals (see corollary \ref{c:21} in \S\,\ref{ss:21}).

\subsubsection{Commuting vector fields and the Cauchy-Riemann equations.}
Let $Z=f(z)\frac{\pa}{\pa z}$ be a germ at $0\in\C$ of holomorphic vector field.
We can split $Z$ into real and imaginary parts, $Z=X+i\,Y$, where $X$ and $Y$ are germs at $0\in\R^2$ of real analytic vector fields. Explicitly, if we write $z=x+i\,y$, $f(z)=u(x,y)+i\,v(x,y)$ and $\frac{\pa}{\pa z}=\frac{\pa}{\pa x}-i\,\frac{\pa}{\pa y}$, then $Z=X+i\,Y$, where $X=u(x,y)\frac{\pa}{\pa x}+v(x,y)\frac{\pa}{\pa y}$ and $Y=v(x,y)\frac{\pa}{\pa x}-u(x,y)\frac{\pa}{\pa y}$.

The Cauchy-Riemann equations imply that $X$ and $Y$ commute: $[X,Y]=0$. Since $X$ and $Y$ are real analytic we can consider their complexification, that we will denote again by $X$ and $Y$.
This classical construction gives examples of non-homogeneous pairs of commuting vector fields on $(\C^2,0)$ with $r(X)=2$.
Up to conjugacy, the vector field $f(z)\frac{\pa}{\pa z}$ can be written as

a. $\frac{\pa}{\pa z}$, if $f(0)\ne0$. In this case, $X$ and $Y$ are two commuting constant vector fields.

b. $\mu\,z\frac{\pa}{\pa z}$, $f(0)=0$ and $f^\pr(0)=\mu\ne0$. In this case, $X$ and $Y$ are two linear vector fields with eigenvalues $\mu=a+i\,b$, $\ov\mu=a-i\,b$ and $\la=b+i\,a$, $\ov\la=b-i\,a$, respectively.

c. $X_{p,\,\la}=\frac{z^{p+1}}{1+\la\,z^p}\frac{\pa}{\pa z}$, if $f$ has multiplicity $p+1\ge2$ at the origin.
In this case, if $\la=0$ then $X$ and $Y$ are homogeneous and have a non-constant meromorphic first integral. If $\la\ne0$ then $X$ and $Y$ have no meromorphic first integral.

\vskip.1in
In fact, the above construction can be generalized in any dimension. Given $Z\in\X_n$ we can write $Z=X+i\,Y$, where $X$ and $Y$ are germs of real analytic vector fields on $(\R^{2n},0)$. Their complexifications, called still $X$ and $Y$, are two commuting vector fields on $(\C^{2n},0)$.
If $Z(0)=0$ then the distribution generated by $X$ and $Y$ has a singular locus of dimension $\ge1$: the set $\{z\in(\C^{2n},0)\,|\,X(z)\wedge Y(z)=0\}$.

For instance, if $n=2$ then the Camacho-Sad theorem on the existence of an analytic separatrix of $Z$ (cf. \cite{cs}), gives a holomorphic separatrix $\g$ for $Z$ through $0\in\C^2$. The germ of surface obtained by complexification of $\g$, considered as real surface on $\R^4\simeq\C^2$, is invariant by both vector fields $X$ and $Y$. This motivates the following problem:

\begin{prob}\label{pr:1}
{\rm Let $X$ and $Y$ be two germs at $0\in\C^n$, $n\ge3$, of commuting vector fields such that $X\wedge Y\not\equiv0$. Does there exists a germ a complex surface through $0\in\C^n$ which is simultaneously $X$ and $Y$ invariant.}
\end{prob}

In the case $n=3$ problem \ref{pr:1} has positive answer \cite{rr}.
However, we would like to note that in general, if $X\wedge Y\not\equiv0$ then there is an analytic set of dimension $\ge1$ invariant by both vector fields, $X$ and $Y$.
This is a consequence of the following result:

\begin{prop}
Let $X\in\X_n$ and assume that the unique proper analytic subset of $(\C^n,0)$ which is $X$-invariant is the origin $\{0\}$.
Then $C(X)=\C.\,X$.
\end{prop}

{\it Proof.}
The hypothesis implies that $0\in Sing(X)$: if not, then the orbit of $0$ by the local flow of $X$ has dimension one and is $X$-invariant. Note also that $0$ is an isolated singularity of $X$.

Assume by contradiction that $C(X)\ne\C.\,X$. In this case, we have two possibilities:
\begin{itemize}
\item[$1^{st}$.] $r(X)=1$ and $d(X)>1$.
\item[$2^{nd}$.] $r(X)>1$.
\end{itemize}
In the first case we have $\I(X)\ne\C$, which is not possible with the hypothesis.
In the second case, there exists $Y\in C(X)$ such that $X\wedge Y\not\equiv0$.
Consider the analytic set
\[
\Si=\{z\in(\C^n,0)\,|\,X(z)\wedge Y(z)\ne0\}\,.
\]
Note that $dim(\Si)\ge1$, because $0\in \Si$.
On the other hand, if $z_o\ne0$ and $z_o\in\Si$ then there are local coordinates $(U,x=(x_1,...,x_n))$ at $z_o$ such that $x(z_o)=0$ and:
\begin{itemize}
\item[(a).] $X|_U=\frac{\pa}{\pa x_1}$, because $X(z_o)\ne0$.
\item[(b).] $Y(x)=\sum_{j=1}^na_j(x_2,...,x_n)\frac{\pa}{\pa x_j}$, because $[X,Y]=0$.
\end{itemize}
Since $z_o\in\Si$ we get
\[
X(0)\wedge Y(0)=0\,\,\implies\,\,a_j(0)=0\,,\,\forall\,j\ge2\,\,\implies\,\,X\wedge Y(x_1,0,...,0)=0\,\,\implies
\]
the analytic subset $A=(x_2=...=x_n=0)\sub U$ is contained in $\Si$ and is $X$-invariant, contradicting the hypothesis.
\qed
\vskip.1in
As an imediate consequence we have:
\begin{cor}
Let $X\in\X_n$ and assume that $d(X)\ge 2$. Then there exists a proper $X$-invariant analytic subset $\Si\sub(\C^n,0)$ with $dim(\Si)\ge1$.
\end{cor}

\section{The case of dimension two}\label{ss:3}
In this section we study the case of dimension two.
A crucial fact that will be used is that if $X\in\X_2$ has a non-constant holomorphic first integral, then there exists $f\in\I(X)$ such that $\I(X)=\C\{f\}$ (see \cite{mm}). The first integral $f$ is called a {\it minimal} first integral of $X$.

Another remark is that when $X\in\X_2$ has an isolated singularity at $0\in\C^2$, $r(X)=1$ and $d(X)\ge2$, then $X$ has a non-constant holomorphic minimal first integral $f$, so that $C(X)=\C\{f\}.\,X$ and $d(X)=\infty$.

In fact, since $r(X)=1$ and $d(X)\ge2$ there exists $Y\in C(X)$ such that $Y\notin\C.\,X$.
But, $r(X)=1$ implies that $Y\wedge X=0$ and so $Y=g.\,X$ where $g\in\O_2\setminus\C$ by Hartogs extension theorem (here we use that $X$ has an isolated singularity at $0\in\C^2$). Note that $X(g)=0$ because
\[
0=[X,Y]=[X,g.\,X]=X(g).\,X\,\,\implies\,\,X(g)=0\,.
\]

\subsection{The case $r(X)=2$ in dimension two.}\label{ss:31}
We have proved in proposition \ref{p:3} that if $X\in\X_2$ and $r(X)=2$ then $X$ is Liouvillian integrable: the foliation $\fa_X$ is also defined by a germ of closed meromorphic 1-form. In other words, the dual form of $X$, $\om:=i_Xdx\wedge dy$, has an integrating factor
$g$:
$d\left(\frac{\om}{g}\right)=0$.
For instance, the saddle-node in the formal normal form is Liouvillian integrable.

In the general case, if the decomposition of the integrating factor $g$ into irreducible factors is
$\Pi_{j=1}^rf_j^{k_j}$ then it is proved in \cite{cema} that:
\begin{equation}\label{eq:om}
\frac{\om}{f}=\sum_{j=1}^r\la_i\,\frac{df_i}{f_i}+d\left(\frac{\var}{f_1^{k_1-1}...f_r^{k_r-1}}\right)\,,
\end{equation}
where
\begin{itemize}
\item[1.] $f_1,...,f_r\in\O_2$ are the poles of $\om$, $\var\in\O_2$ and $f_j$ does not divides $\var$, $1\le j\le r$.
\item[2.] $\la_j\in\C$ is the residue of $\frac{\om}{f}$ along $f_j$, $1\le j\le r$.
\item[3.] $k_j\in\N$ is the multiplicity of $\frac{\om}{f}$ at the pole $f_j$, $1\le j\le r$.
\end{itemize}
The foliation $\fa_X$ is said to be {\it logarithmic} if $\var/f_1^{k_1-1}...f_r^{k_r-1}$ is holomorphic, or equivalently $k_1=...=k_r=1$.
\vskip.1in
\begin{rem}\label{r:31}
{\rm Each irreducible component of the pole $f$ of $\frac{\om}{f}$ defines a separatrix $\Ga_j:=(f_j=0)$ of $\fa_X$.

Another observation is that the "multivalued" function $\sum_j\la_j\,log(f_j)+\var/f^{k_1-1}...f_r^{k_r-1}$ is a first integral of $X$.}
\end{rem}
\vskip.1in
Suppose now that $X_1$ and $X_2$ are two holomorphic germs on $(\C^2,0)$ of commuting vector fields such that $X_1\wedge X_2\not\equiv0$.
\begin{claim}\label{cl:31}
There are unique closed meromorphic 1-forms $\a_1$ and $\a_2$ such that $\a_i(X_j)=\d_{ij}$, $i,j=1,2$, where $\d_{ij}=1$ if $i=j$ and $0$ otherwise.
\end{claim}

{\it Proof.}
As we have seen in the proof of proposition \ref{p:3}, if we set $X_1\wedge X_2=g\frac{\pa}{\pa z_1}\wedge\frac{\pa}{\pa z_2}$, $\om_1=i_{X_2}dz_1\wedge dz_2$ and $\om_2=-i_{X_1}dz_1\wedge dz_2$, then the forms $\a_1:=\frac{\om_1}{g}$ and $\a_2:=\frac{\om_2}{g}$ are closed.
As the reader can check, we have also $\a_i(X_j)=\d_{ij}$, $1\le i,j\le2$.
\qed
\vskip.1in
As a consequence, we have the following:
\begin{prop}\label{p:4}
Given $X\in\X_2\setminus\{0\}$ we have three possibilities:
\begin{itemize}
\item[(a).] $r(X)=d(X)=1$. In this case, $C(X)=\C.\,X$.
\item[(b).] $d(X)\ge2$ and $r(X)=1$. In this case, $X$ has a non-constant meromorphic first integral. 
\item[(c).] $r(X)=2$. In this case $X$ is Liouvillian integrable.
\end{itemize}
\end{prop}

\begin{cor}\label{c:31}
Let $X\in\X_2$ be a saddle-node. Then $X$ is holomorphically normalisable if, and only if, $r(X)=2$.
\end{cor}

\begin{cor}
Let $X\in\X_2$ with $DX(0)=x\frac{\pa}{\pa x}+\la\,y\frac{\pa}{\pa y}$, with $\la\in\C^*\setminus\Q$.
Then $X$ is linearizable if, and only if, $r(X)=2$.
\end{cor}

Another interesting consequence is the following:

\begin{cor}
Let $X\in\X_2$ with $r(X)=2$. Let $\Pi\colon(M,E)\to(\C^2,0)$ be the Seidenberg resolution of singularities of $X$ and $\wt\fa_X$ be the strict transform of the foliation defined by $\Pi^*(X)$ (cf. \cite{se} and \S\,\ref{ss:32}). Then the holonomy of $\wt\fa_X$ in any non-dicritical irreducible component of $E$ is abelian.
\end{cor}

\subsection{Germs of vector fields with just one irreducible separatrix.}\label{ss:32}

In the case of dimension two, an important class of foliations on $(\C^2,0)$ are the so-called "generalized curves" defined in \cite{cls}. These foliations have no saddle nodes in their reduction of singularities. According to Seidenberg`s resolution theorem \cite{se}, given a germ at $0\in\C^2$ of foliation $\fa$ there exists a blowing-up procedure $\Pi\colon(M,D)\to(\C^2,0)$ such that the strict transform $\Pi^*(\fa)$ of $\fa$ has only {\it reduced} singularities on the exceptional divisor $D$.
We say that a singularity $p\in D$ of $\Pi^*(\fa)$ is reduced if the germ at $p$ of the foliation is represented by a vector field $Y$ such that its linear part $DY(p)$ at $p$ has eigenvalues $\la_1,\la_2$, where:
\begin{itemize}
\item{} either $\la_1=0$ and $\la_2\ne0$ (or vice-versa). In this case, the singularity is a saddle-node,
\item{} or $\la_1.\la_2\ne0$ and $\la_2/\la_1\notin\Q_+$.
\end{itemize}
We will say that $X\in\X_2$ is a generalized curve (briefly G.C) if the associated foliation has no saddle-node in its resolution of singularities.
In this case, we have the following result:
\begin{thm}\label{t:1}
Let $X\in\X_2$ be a G.C vector field with an isolated singularity at $0\in\C^2$ and with just one separatrix. Then $r(X)=1$. In particular,
$C(X)=\I(X).\,X$. 
\end{thm}

When a germ $X\in\X_2$ has a holomorphic first integral $f$, whoose decomposition into irreducible factors is $f=\Pi_{j=1}^rf_j^{k_j}$, then the separatrices of $X$ through $0\in\C^2$ are the curves
$(f_j=0)$, $1\le j\le r$. In this case, $X$ is non-dicritical.
In particular, if $X$ has an irreducible holomorphic first integral, say $f$ , then $(f=0)$ is the unique irreducible separatrix of $X$. In this case, we have the following:

\begin{cor}\label{c:32}
Let $X\in\X_2$ and assume that $0$ is an isolated singularity of $X$ and that $X$ has an irreducible first integral $f\in\O_2$.
Then $r(X)=1$ and $C(X)=\C\{f\}.\,X$. Moreover, there exists a local coordinate system $(x,y)$ around $0\in\C^2$ such that $X=H_f$, where $H_f=f_y\frac{\pa}{\pa x}-f_x\frac{\pa}{\pa y}$ is the Hamiltonian of $f$.
\end{cor}

Before proving theorem \ref{t:1} we will see two examples in which the conclusion of theorem \ref{t:1} is not true.

\begin{ex}\label{ex:1}
{\rm If $X$ is non singular, $X(0)\ne0$, then in some coordinate system $(x,y)$ around $0$ we have $X=\frac{\pa}{\pa x}$. On the other hand, it is not difficult to see that
\[
C\left(\frac{\pa}{\pa x}\right)=\C\{y\}\frac{\pa}{\pa x}+\C\{y\}\frac{\pa}{\pa y}\,.
\]}
\end{ex}

\begin{ex}\label{ex:2}
{\rm If $f$ is not irreducible then the theorem is not true in general. For instance, if
$X=x\frac{\pa}{\pa x}-y\frac{\pa}{\pa y}=H_{x.y}$ then $x.y$ is a minimal first integral of $X$: $\I(X)=\C\{x.y\}$, but 
\[
C(X)=\C\{x.y\}.\,x\frac{\pa}{\pa x}\oplus\C\{x.y\}.\,y\frac{\pa}{\pa y}\,,
\]
which is not of the form $\C\{f\}.\,X$.}
\end{ex}

{\it Proof of theorem \ref{t:1}.}
Observe first that the irreducibility of $f$ implies that if $\I(X)\ne\C$ then, modulo an unit, we can assume that $f$ is a minimal first integral of $X$
(see \cite{mm}). Since $f\in\I(X)$ the curve $\Ga:=(f=0)$ is $X$-invariant (it is a separatrix of $X$).  In fact, we have the following:
\begin{lemma}\label{l:31}
If $Y\in C(X)$ then $\Ga$ is $Y$-invariant. Moreover, there exists $\la\in\C$ such that
$Y|_\Ga=\la.\,X|_\Ga$.
\end{lemma}
\vskip.1in
{\it Proof.} 
We will assume $Y\not\equiv0$.
In some coordinate system $(x,y)$ around $0\in\C^2$ we can write $X=X_1\frac{\pa}{\pa x}+X_2\frac{\pa}{\pa y}$ and
\[
X\wedge Y=g.\,\nu\,,\,\,\text{where}\,\,\nu=\frac{\pa}{\pa x}\wedge\frac{\pa}{\pa y}\,\,\text{and}\,\,g\in\O_2
\]
and $g(0)=0$.
\vskip.1in
$1^{st}$ case. $g\equiv0$. In this case, we have $Y=\phi.\,X$, where $\phi\in\O_2$, because $X$ has an isolated singularity at $0$.
Since $Y\in C(X)$ we get $Y(\phi)=0$ and $\phi\in\I(X)$.
If $\phi\in\C^*$ we are done. If $\phi$ is non-constant then $\I(X)\ne\C$ and we can assume that $f\in\I(X)$.

From $[X,Y]=0$ we get
\[
X(Y(f))=Y(X(f))=0\,\,\implies\,\,Y(f)\in\I(X)\,\implies\,f\,|\,Y(f)\,\,\implies
\]
$\Ga$ is $Y$-invariant.
\vskip.1in
$2^{nd}$ case. $g\not\equiv0$.
If $L$ denotes the Lie derivative then
\[
0=L_X\,(X\wedge Y)=[X(g)+\nabla(X).\,g]\,\nu\,,\,\text{where}\,\nabla(X)=\frac{\pa X_1}{\pa x}+\frac{\pa X_2}{\pa y}\,,
\]
which implies $g\,|\,X(g)$ and so $g$ is $X$-invariant and each irreducible component of $g$ is a separatrix of $X$.
Since $\Ga$ is the unique separatrix of $X$ we must have $g=u.\,f^k$, where $u\in\O_2^*$ and $k\ge1$.
Analogously,
\[
0=L_Y\,(X\wedge Y)=[Y(g)+\nabla (Y).\,g]\,\nu\,\,\implies
\]
$\Ga$ is $Y$-invariant.
\vskip.1in
Let $\g(t)=(x(t),y(t))$ be a Puiseux`s parametrization of $\Ga$. Since $\Ga$ is $X$ and $Y$-invariant we can write $\g^*(X)=\phi(t)\frac{\pa}{\pa t}$ 
and $\g^*(Y)=\var(t)\frac{\pa}{\pa t}$. Note that $\phi(0)=0$, but $\phi\not\equiv0$, for otherwise $0$ would not be an isolated singularity of $X$. On the other hand, $[X,Y]=0$ implies
\[
\left[\phi(t)\,\frac{\pa}{\pa t},\var(t)\,\frac{\pa}{\pa t}\right]=0\,\implies\,\phi(t).\,\var^\pr(t)-\var(t).\,\phi^\pr(t)=0\,\implies
\]
\[
\var(t)=\la.\,\phi(t)\,,\,\la\in\C\,,
\]
which proves lemma \ref{l:31}.
\qed
\vskip.1in
{\it End of the proof of theorem \ref{t:1}.}
Let $Z_1=Y-\la.\,X$. 
From $Z_1|_\Ga=0$ we get $Z_1=f^{k_1}.\,Y_1$, where $Y_1\in\X_2$ and $k_1\ge1$.
We have two possibilities:
\vskip.1in
{\it $1^{st}$ case:} $\I(X)\ne\C$. In this case, we can assume that $\I(X)=\C\{f\}$.
Since $Z_1\in C(X)$ and $f\in \I(X)$ we get $Y_1\in C(X)$. From lemma \ref{l:31} we get again
$Y_1=\la_1.\,X+f^{k_2}.\,Y_2$, where $\la_1\in\C$, $Y_2\in\X_2$ and $k_2\ge1$, so that
\[
Y=\left(\la+\la_1.\,f^{k_1}\right).\,X+f^{k_1+k_2}.\,Y_2\,,
\]
where $[X,Y_2]=0$. Using the above argument inductively we get a formal series $g(z)=\sum_{j=0}^\infty\la_j\,z^j$ such that $Y=g(f)\,.X$. Since $X$ and $Y$ are holomorphic, it is clear that the series $g$ is convergent and so $Y\in \C\{f\}.\,X$ and $r(X)=1$.
\vskip.1in

{\it $2^{nd}$ case:} $\I(X)=\C$. We have two sub-cases:

2.a. $X\wedge Y\equiv0$. In this case, we have $Y=h.\,X$ where $h\in\O_2$, because $X$ has an isolated singularity at $0\in\C^2$.
Here $h\in\I(X)$ and since $\I(X)=\C$ we get $h\in\C$ and so $h=\la$ and $Y=\la.\,X\in\C.\,X$.

2.b. $X\wedge Y\not\equiv0$. We will see that $\I(X)\ne\C$, a contradiction.

In fact, let $Y\in C(X)$ be such that $X\wedge Y\ne0$. 
We assert that the set $\Ga:=\{z\in(\C^2,0)\,|\,X(z)\wedge Y(z)=0\}$ is $X$-invariant:

The set $\Ga$ is an analytic curve through $0\in\C^2$: if we fix local coordinates $(x,y)$ around $0$, then $X\wedge Y=g.\,\mu$, where $\mu:=\frac{\pa}{\pa x}\wedge\frac{\pa}{\pa y}$, and $\Ga=(g=0)$. On the other hand,
\[
L_XX\wedge Y=[X,X]\wedge Y+X\wedge[X,Y]=0\,\implies
\]
\[
0=L_X(g\,\mu)=X(g)\,\mu+g.\,L_X\mu=X(g)\,\mu-\nabla X\,g\,\mu\,\implies\,X(g)=\nabla X.\,g\,,
\]
where $\nabla X=\frac{\pa X(x)}{\pa x}+\frac{\pa X(y)}{\pa y}$. Therefore, $\Ga$ is $X$-invariant and is the unique separatrix of $X$. In particular, $\Ga=(f=0)$ and after multiplying $g$ by an unit we can assume $g=f^k$. 
Let $\a$ and $\be$ be the meromorphic closed 1-forms with $\a(X)=\be(Y)=0$ and $\a(Y)=\be(X)=1$. As we have seen before $\a=\frac{1}{f^k}\,i_X\nu$ and $\be=\frac{1}{f^k}\,i_Y\nu$, $\nu=dx\wedge dy$.
In particular, we can write \cite{cema}
\[
\a=\la\,\frac{df}{f}+d\left(\frac{H}{f^{k-1}}\right)\,,
\]
where $\la\in\C$, $H\in\O_2$, and, either $k=1$ and $\la\ne0$, or $k>1$.
When $k>1$ we will assume also that $f$ does not divides $H$.

If $\la\ne0$ and $k=1$ then $\wt{f}:=f.\,e^{\la^{-1}H}$ is an irreducible holomorphic first integral of $X$, because
\[
\frac{d\wt{f}}{\wt{f}}=\frac{df}{f}+\la^{-1}dH=\frac{1}{\la}\,\a\,.
\]

If $k>1$ then $H\not\equiv0$ and we can write $H=h_0+h$, where $h_0=H(0)$. We have four possilities:
\begin{itemize}
\item[$1^{st}$.] $h_0=0$ and $\la=0$. In this case, $X$ is dicritical, because $f$ does not divides
$H=h$ and all curves of the form $(f^{k-1}-c\,h=0)$, $c\in\C$, are $X$-invariant, so that $X$ has more than one separatrix.
\item[$2^{nd}$.] $h_0\ne0$ and $\la=0$. In this case, there exists $\phi\in\O_2^*$ such that $\phi^{k-1}=H$ and $\wt{f}:=f/\phi$ is an irreducible holomorphic first integral of $X$.
\item[$3^{rd}$.] $h_0\ne0$ and $\la\ne0$. In this case, we have
\begin{equation}\label{eq:om1}
f^k.\,\a=\la\,f^{k-1}\,df-(k-1)(h_0+h)df+f\,dh
\end{equation}
\end{itemize}
We see that $f^k\,\a$ is the pull-back by the morphism $(x,y)\mapsto(f,h)=(u,v)$ of the 1-form
\[
\wt\a:=\la\,u^{k-1}du+u\,dv-(k-1)(h_0+v)\,du\,.
\]
Since $\wt\a(0)=-(k-1)\,h_0\,du\ne0$, by Frobenius theorem $\wt\a$ has a holomorphic first integral
$F(u,v)$ of the form $(u,v)\mapsto F(u,v)=u+h.o.t.$, so that $X$ has a first integral of the form
$F(g,h)=f+...$. Since $(f=0)$ is a separatrix of $X$, we must have $F(f,h)=u.\,f:=\wt{f}$, 
$u\in\O_2^*$, so that
$X$ has an irreducible holomorphic  first integral.

$4^{th}$. $h_0=0$ and $\la\ne0$. In this case (\ref{eq:om1}) can be written as
\[
f^k.\,\a=\la\,f^{k-1}\,df-(k-1)\,h\,df+f\,dh
\] 
and $f^{k-1}.\,\a$ is the pull-back of the 1-form
\[
(\la\,u^{k-1}-(k-1)v)\,du+u\,dv
\]
by the map $(x,y)\mapsto(f,h)=(u,v)$. The dual vector field of the above 1-form is
\[
Y=u\,\frac{\pa}{\pa u}+\left((k-1)\,v-\la\,u^{k-1}\right)\,\frac{\pa}{\pa v}
\]
which is in the Poincar\'e-Dulac normal form. If $\la\ne0$ then this vector field has a saddle-node in its redution of singularities \cite{mm}.
This implies that $X$, the original vector field, has the same property and so is not G.C.
Hence, in all possible cases, $X$ has an irreducible holomorphic first integral.
\qed
\vskip.1in
{\it Proof of corollary \ref{c:32}.}
Since $X$ has a holomorphic first integral $f$ it is G.C (cf. \cite{cls}). Since $f$ is irreducible 
$X$ has an unique separatrix.
Therefore, by theorem \ref{t:1} we have $C(X)=\C\{f\}.\,X$.
It remains to prove that in some coordinate system aroud $0\in\C^2$ we have $X=H_f$. 

Fix a local coordinate system $(u,v)$. Since $X(f)=0$ and $0\in\C^2$ is an isolated singularity of $X$ and of $df$, there exists an unity $\phi\in\O_2^*$ such that
$\phi.\,i_X\,du\wedge dv=df$. If we consider the change of variables $x=u$ and $y=\var.\,v$, where $\var\in\O_2^*$, then we get
\[
dx\wedge dy=(v.\,\var_v+\var)\,du\wedge dv\,.
\]
The p.d.e. $v.\,\var_v+\var=\phi$ has a solution $\var\in\O_2^*$. Hence in the new coordinate system we have
\[
i_X\,dx\wedge dy=df\,\,\implies\,\,X=H_f\,\,.\qed
\]

\begin{rem}
{\rm Given $f,g\in\O_2\setminus\{0\}$, with $f(0)=g(0)=0$, it is easy to check that $[H_f,H_g]=0$ if, and only if, the jacobian determinant of $det\,J(f,g)$ is a constant $\la$. When $\la\ne0$ then the map $(f,g)\colon(\C^2,0)\to(\C^2,0)$ is a germ of biholomorphism. In particular, $f$ is a submersion and $H_f(0)\ne0$.

If $\la=0$ then, after \cite{mm}, we can write $f=\phi(h)$ and $g=\psi(h)$, where $\phi,\psi\in\O_1$ and $h$ is a primitive first integral of $H_f$. We are essentially in the situation of theorem
\ref{t:1}: $H_f, H_g\in C(H_h)$ and $r(H_f)=1$.

In contrast, in dimension $n\ge3$ we have examples of $f\in\O_n$ irreducible, but with $r(H_f)>1$.
This type of example can be constructed as follows: let $\C^n=E_1\oplus...\oplus E_s$ be a linear decomposition of $\C^n$, $X_i\in\X(E_i,0)$, $i=1,...,s$, and set $X=X_1+...+X_s\in\X_n$.
Then it is clear that
\[
C(X)\sup\bigoplus_{j=1}^s\I(X_j).\,X_j\,,
\]
where $\I(X_j)\sub\O(E_j,0)$ is the subring of first integrals of $X_j$, $1\le j\le s$.

For instance, let $\C^{2n}=\C^2\oplus...\oplus\C^2$ with the sympletic form $\Om=dx_1\wedge dy_1+...+dx_n\wedge dy_n$.
For each $j=1,...,n$ let $f_j=f_j(x_j,y_j)$ and $f=f_1+...+f_n$, so that the associated hamiltonian is given by $df=i_{H_f}\Om$,
\[
H_f=\sum_{j=1}^n\frac{\pa f_j}{\pa y_j}\frac{\pa}{\pa x_j}-\frac{\pa f_j}{\pa x_j}\frac{\pa}{\pa y_j}
\]
The ring $\r:=\C\{f_1,...,f_n\}$ is a subring of $\I(H_f)$ and $C(H_f)$ contains
\[
\bigoplus_{j=1}^n\r.\,\left(\frac{\pa f_j}{\pa y_j}\frac{\pa}{\pa x_j}-\frac{\pa f_j}{\pa x_j}\frac{\pa}{\pa y_j}\right)
\]
In fact, given $h_1,h_2\in\r$ then $[H_{h_1},H_{h_2}]=0$.}
\end{rem}

\subsection{Generalized curves.}\label{ss:24}
In this section we will study commuting G.C vector fields of a particular type.

\begin{definition}\label{def:1}
{\rm We say that $X\in\X_2$ has an {\it isolated non-rational} singularity at $0\in\C^2$ if its linear part $DX(0)$ has non-vanishing eigenvalues, $\la_1,\la_2$, with $\la_2/\la_1\notin\Q$.
When the quotient $\la_2/\la_1\in\Q$ we will say that $0$ is an isolated rational singularity.

We say that $Y\in\X_2$ has a {\it non-rational} singularity at $0\in\C^2$ if $Y=g.\,X$, where $g\in\O_2$ and $X$ has an isolated non-rational singularity at $0$.

If $X\in\X_2$ has a non-reduced singularity at $0\in\C^2$ then we say that it is of {\it non-rational} type if it has a non-rational singularity in its reduction of singularities.}
\end{definition}

Our main result in this section is the following:

\begin{thm}\label{t:2}
Let $X\in\X_2$ be a G.C (non-dicritical) with an isolated and non-rational singularity at
$0\in\C^2$. If $DX(0)$ is nilpotent then $C(X)=\C.\,X$.
\end{thm}

{\it Proof.}
We need a lemma.
\begin{lemma}\label{l:32}
Let $Z\in\X_2$ with a singularity at $0\in\C^2$ and $\Pi\colon(\wt\C^2,D)\to(\C^2,0)$ be a blowing-up with exceptional divisor $D\simeq\p^1$. If $DZ(0)=0$ then:
\begin{itemize}
\item[(a).] $\Pi^*(Z)|_D\equiv0$.
\item[(b).] For any singularity $p\in D$ of the strict transform of $\Pi^*(Z)$ then $D\Pi^*(Z)(p)=0$.
\end{itemize}
\end{lemma}

{\it Proof.} The proof relies in the Seidenberg's reduction process of $Z$ \cite{se}. Write the Taylor series of $Z$ as
\[
Z=Z_\nu+Z_{\nu+1}+...=\sum_{j\ge\nu}Z_j\,\,,
\]
where $Z_j$ is homogeneous of degre $j$ and $Z_\nu$ is the first non-zero jet of $Z$. Since $DZ(0)=0$ we have $\mu(Z,0)=\nu\ge2$.
Locally, in suitable coordinates, $\Pi$ is a quadratic map $\Pi(x,t)=(x,t.\,x)$ and $\Pi^{-1}(0)=(x=0)$. If $Z=A(x,y)\frac{\pa}{\pa x}+B(x,y)\frac{\pa}{\pa y}$ then
\[
\Pi^*(Z)=A(x,t.\,x)\frac{\pa}{\pa x}+\frac{B(x,t.\,x)-t.\,A(x,t.\,x)}{x}\,\frac{\pa}{\pa t}\,.
\]
If $Z=\sum_{j\ge\nu}Z_j$, as before, where $Z_j=A_j(x,y)\frac{\pa}{\pa x}+B_j(x,y)\frac{\pa}{\pa y}$ then we get $\pi^*(Z)=\sum_{j\ge\nu}\Pi^*(Z_j)$, where
\[
\Pi^*(Z_j)(x,t)=x^{j-1}.\,\left[x.\,A_j(1,t)\frac{\pa}{\pa x}+(B_j(1,t)-t.\,A_j(1,t))\frac{\pa}{\pa t}\right]:=x^{j-1}.\,\wt{Z}_j(x,t)\,.
\]
From the above formula, we have two possibilities:
\begin{itemize}
\item[$1^{st}$.] $B_\nu(1,t)-t.\,A_\nu(1,t)\not\equiv0$. We see that $\Pi^*(Z)=x^{\nu-1}.\,\wt{Z}$, where $\wt{Z}$ has isolated singularities on the exceptional divisor $D$, which is $\wt{Z}$-invariant. This is case of a non-dicritical blow-up. The foliation on $\wt\C^2$ induced by $\wt{Z}$ is called the strict transform of $Z$ and denoted as $\wt\fa_Z$.
Observe also that
\begin{itemize}
\item[1.1.] $Sing(\wt\fa_Z)\cap D\ne\emp$.
\item[1.2.] If $p\in Sing(\wt\fa_Z)\cap D$ then the algebraic multiplicity of $\Pi^*(Z)$ at $p$, $\mu(\Pi^*(Z),p)\ge\nu$. 
\end{itemize}
\item[$2^{nd}$.] $B_\nu(1,t)-t.\,A_\nu(1,t)\equiv0$. We see that $\Pi^*(Z)=x^\nu.\,\wt{Z}$, where $\wt{Z}$ has isolated singularities on $D$, which is not $\wt{Z}$-invariant. This is the case of a dicritical blow-up. Note that this happens if, and only if, $X_\nu$ is colinear with the radial vector field $x\frac{\pa}{\pa x}+y\frac{\pa}{\pa y}$.
\end{itemize}
This proves lemma \ref{l:32}.
\qed

Let us continue the proof of theorem \ref{t:2}.
In the non-dicritical hypothesis of theorem \ref{t:2} we are excluding the $2^{nd}$ possibility in any step of the reduction process for $X$.
Let $\Pi\colon(M,D)\to(\C^2,0)$ be the reduction of the vector field $X$, where $D:=\Pi^{-1}(0)$ is the exceptional divisor. Let $D=\bigcup_{j=1}^rD_j$ be the decomposition of $D$ into irreducible components. Denote by $\wt\fa_X$ the foliation induced by $\Pi^*(X)$.
We will assume that $D_1$ is strict transform by $\Pi$ of the divisor obtained at the first step of the resolution. 

\begin{rem}\label{re:33}
{\rm When $DX(0)\ne0$ is nilpotent we can assume $DX(0)=y\frac{\pa}{\pa x}$. In this case, in the first blow-up $\Pi_1\colon(\wt\C^2,D_1)\to(\C^2,0)$ we find just one singularity on the divisor $D_1$, the point $p\in D_1$ corresponding to the direction $y=0$. The vector field $\Pi_1^*(X)$ is not identically zero on $D_1$, but $\mu(\Pi_1^*(X),p)=2$ (see \cite{mm}). The point $p$ is not a simple singularity for the foliation defined by $\Pi_1^*(X)$ and so we have to blow-up more times in the resolution process. If $X$ is G.C and at each step of the resolution we blow-up only at singularities of the strict transform then lemma \ref{l:32} implies that $\Pi^*(X)|_{D_j}\equiv0$ for all $j\ge2$. As a consequence, we have the following:
\begin{cor}
Let $D_k$ be an irreducible component of $D$ where there is a non-rational singularity of $\wt\fa_X$. Then $\Pi^*(X)|_{D_k}\equiv0$. 
\end{cor}}
\end{rem}

In the above resolution if $i\ne j$ then, either $D_i\cap D_j=\emp$, or $D_i$ cuts $D_j$ transversely in just one point.
Denote as $\wt\fa_X$ the strict transform of the foliation induced by the vector field $\Pi^*(X)$.
As a consequence of the above computation, we can conclude that:
\begin{itemize}
\item[(i).] The strict transform $\wt{\fa}_X$ has singularities in all components $D_j$ of the exceptional divisor (see \cite{br} and \cite{mm}). This follows from the hypothesis that $X$ is non-dicritical.
\item[(ii).] The vector field $\Pi^*(X)$, which is holomorphic on $M$, vanishes identically along all components $D_j$ of exceptional divisor, except perhaps at $D_1$.
\end{itemize}

\begin{claim}\label{cl:32}
Let $p\in D$ be a non-rational singularity of $\wt\fa_X$. Then:
\begin{itemize}
\item[(a).] There are local coordinates around $p$, $(x,y)\colon(M,p)\to(\C^2,0)$, such that the first non-zero jet of $\Pi^*(X)$ is of the form
$x^p\,y^q\,\left(x\frac{\pa}{\pa x}+\la\,y\frac{\pa}{\pa y}\right)$, where
$\la\notin\Q$ and $p+q\ge1$.
\item[(b).] The germ of $\Pi^*(X)$ at $p$ is formally conjugated to $x^p\,y^q\,\left(x\frac{\pa}{\pa x}+\la\,y\frac{\pa}{\pa y}\right)$.
\end{itemize}
\end{claim}

{\it Proof.} 
First of all, we have two possibilities for $p\in D$:
\begin{itemize}
\item[(1).] $p\in D_j$ and $p\notin D_i$ if $i\ne j$.
\item[(2).] $p\in D_i\cap D_j$, with $i\ne j$, and $p\notin D_s$, $s\ne i,j$.
\end{itemize}

Note that there are local coordinates around $p$, $(x,y)\colon(U,p)\to(\C^2,0)$, in which in the first case we have $D\cap U=D_j\cap U=(x=0)$ and in the second case we have $D\cap U=(D_i\cup D_j)\cap U=(x.\,y=0)$. 

It follows from (ii) above that $\Pi^*(X)$ vanishes in a certain order $\ge1$ along each $D_i\sub D$, so that, in both cases we can write $\Pi^*(X)=x^p\,y^q\,Z$, where $p+q\ge1$ and $Z$ has an isolated singularity at $0\in\C^2$.
In case (1) we have $p\ge1$ and $q=0$, whereas in case (2) we have $p,q\ge1$.
Since $p$ is a non-rational singularity of $\wt\fa_X$, then $det(DZ(0))\ne0$ and the eigenvalues $\la_1,\la_2$ of $DZ(0)$ satisfy $\la_2/\la_1\notin\Q$. In particular, $Z$ has two smooth and transverse separatrices through $p$ (see \cite{mm}). In case (2) necessarily these separatrices are contained in $D_i\cup D_j=(x.\,y=0)$, whereas in case (1), one of the separatrices is contained in the exceptional divisor $(x=0)$ and the other, after a holomorphic change of variables, we can assume that is $(y=0)$.
In both cases, after dividing $Z$ by $\la_1$, we can assume that the first jet of $Z$ at $0\in\C^2$ is of the form $x\frac{\pa}{\pa x}+\la\,y\frac{\pa}{\pa y}$, $\la=\la_2/\la_1$.
In this case, by Poincar\'e's linearization theorem \cite{ma}, $Z$ is formally conjugated to its linear part $x\frac{\pa}{\pa x}+\la\,y\frac{\pa}{\pa y}$ by a formal $\Phi\in\wh{Diff}(\C^n,0)$ such that $\Phi(u,v)=(u\,\Phi_1(u,v),v\,\Phi_2(u,v))=(x,y)$:
$\Phi^*(Z)=u\frac{\pa}{\pa u}+\la\,v\frac{\pa}{\pa v}:=\wh{Z}$.
In this case, we must have
\[
\Phi^*(\Pi^*(X))=\phi(u,v).\,u^p\,v^q\,\left(u\frac{\pa}{\pa u}+\la\,v\frac{\pa}{\pa v}\right)=\phi.\,u^p\,v^q\,\wh{Z}\,,
\]
where $\phi\in\wh\O_2^*$, and $\phi(0)=1$.
Consider now a formal diffeomorphism of the form $\Psi(x,y)=(e^\a.\,x,e^{\la\,\a}.\,y)=(u,v)$, where $\a(0)=0$.
With a straightforward computation we have
\[
\wh{Z}(x)=\wh{Z}\left(e^{-\a}\,u\right)=\left(1-\wh{Z}(\a)\right)x\,\,\text{and}\,\,
\wh{Z}(y)=\wh{Z}\left(e^{-\la\,\a}\,v\right)=\la\left(1-\wh{Z}(\a)\right)y\,,
\]
so that
\[
\Psi^*\left(\phi.\,u^p\,v^q\,\wh{Z}\right)=\phi.\,e^{(p+\la\,q)\a}\left(1-\wh{Z}(\a)\right)\left(x\frac{\pa}{\pa x}+\la\,y\frac{\pa}{\pa y}\right)\,x^p\,y^q\,.
\]
This reduces the proof of the claim to find a solution of the differential equation
\[
\wh{Z}(\a)=u\,\frac{\pa\a}{\pa u}+\la\,v\,\frac{\pa\a}{\pa v}=1-\phi^{-1}.\,e^{-(p+\la\,q)\a}\,.
\]
If we set $w=e^{(p+\la\,q)\a}$ then the above differential equation becomes
\begin{equation}\label{eq:3}
\wh{Z}(w)-(p+\la\,q)w=-(p+\la\,q)\phi^{-1}\,.
\end{equation}
Since $\la\notin\Q$ the linear operator
\[
w\in\wh\O_2\mapsto\,\wh{Z}(w)-(p+\la\,q)w\,\in\wh\O_2
\]
is surjective, as the reader can check. Hence, the differential equation (\ref{eq:3}) has a formal solution. This proves the claim.
\qed
\vskip.1in
Let us suppose by contradiction that $C(X)\ne\C.\,X$ and let $Y\in C(X)\setminus\C.\,X$. We have two possibilities:
\vskip.1in
$1^{st}$. $Y\wedge X\equiv0$. In this case, since $X$ has an isolated singularity at $0$, we must have $Y=f.\,X$, where $f\in\O_2$ is a non-constant first integral of $X$. But, when $X$ has a non-constant first integral all singularities of the strict transform $\wt\fa_X$ are rational (see \cite{mm}).
\vskip.1in

$2^{nd}$. $Y\wedge X\not\equiv0$. Set $Y\wedge X=f.\,\frac{\pa}{\pa x}\wedge\frac{\pa}{\pa y}$,
$f\in\O_2$, $f(0)=0$. As we have seen in the proof of theorem \ref{t:1} in \S\,\ref{ss:32}, the curve $(f=0)$ is $X$ and $Y$-invariant: if $f=\Pi_{j=1}^rf_j^{k_j}$ is the decomposition of $f$ into irreducible factors, then the curves $\Ga_j=(f_j=0)$, $1\le j\le r$, are all simultaneously
$X$ and $Y$-invariant.
Moreover, if $\a=\frac{1}{f}\,i_Ydx\wedge dy$ and $\be=-\frac{1}{f}\,i_Xdx\wedge dy$ then
$\a(X)=\be(Y)=1$, $\a(Y)=\be(X)=0$ and $d\a=d\be=0$. The curve $\Ga:=(f_1...f_r=0)$ is called the reduced separatrix of $X$.
We denote $\Ga_j=(f_j=0)$, $1\le j\le r$.
\vskip.1in
Let $\Pi\colon(M,D)\to(\C^2,0)$ be the minimal reduction of singularities of $X$. Denote as $\wt\fa_X$ and $\wt\fa_Y$ the foliations given by the strict transforms of $\Pi^*(X)$ and $\Pi^*(Y)$, respectively.

\begin{claim}\label{cl:33}
In the above situation we have
\begin{itemize}
\item[(a).] $\mu(Y,0)\ge\mu(X,0)$.
\item[(b).] $Sing(\wt\fa_X)\cap D\sub Sing(\wt\fa_Y)\cap D$. 
\item[(c).] For any $p\in Sing(\wt\fa_X)\cap D$ we have  $\mu(\Pi^*(X),p)\le\mu(\Pi^*(Y),p)$.
\end{itemize}
\end{claim}

{\it Proof.} Since $X$ is G.C, a result of \cite{cls} says that the reduction of $X$ coincides with the reduction of the curve $\Ga$, which in fact, coincides with the reduction of singularities of the foliation given by $d(f_1...f_r)=0$.
Note that $X$ has no other separatrices than that defined by the $f_{j's}$.
Moreover, in \cite{cls} is also proved that:
\begin{itemize}
\item[(i).] $\mu(X,0)=\mu(d(f_1...f_r),0)=\sum_{j=1}^r\mu(f_j,0)-1$. 
\item[(ii).] If $Z$ is any vector field such that $\Ga_1,...,\Ga_r$ are separatrices of $Z$ then $\mu(Z,0)\ge\mu(X,0)$.
\item[(iii).] If $Z$ is as in (ii) and $D_j$ is any irreducible component of the exceptional divisor $D$, then the order of annulment of $\Pi^*(Z)$ along $D_j$ is $\ge$ the order of annulment of $\Pi^*(X)$ along $D_j$. In other words, if in local coordinates $(x,y)$ near some $z\in D_j$ we have $D_j=(x=0)$, $\Pi^*(X)=x^{k}\,\wt{X}$ and $\Pi^*{Z}=x^\ell\,\wt{Z}$, where $\wt{X}$ and $\wt{Z}$ have isolated singularities, then $k\le\ell$. We will denote $\ell:=\mu(\Pi^*(Z),D_j)$. 
\end{itemize}
In particular, (ii) $\implies$ (a). Let us prove (b) and (c).
Since the reduction of singularities of $X$ and $d(f_1...f_r)=0$ coincide, if $p\in Sing(\wt\fa_X)\cap D$ then:
\begin{itemize}
\item[I.] either $p$ is a corner of $D$: $p=D_i\cap D_j$, $i\ne j$,
\item[II.] or $p$ corresponds to the intersection of the strict transform of some of the curves $\Ga_j$ with some irreducible component $D_i$ of
$D$. In this case, $p\notin D_\ell$ if $\ell\ne i$.
\end{itemize}
In any case, I or II, $p$ must be also a singularity of the strict transform $\wt\fa_Y$, because all the curves $\Ga_j$, $1\le j\le r$, are $Y$-invariant.
Since all singularities of $\wt\fa_X$ are reduced, we have $\mu(\wt\fa_X,p)=1\le \mu(\wt\fa_Y,p)$. 
Finally, in case I, $p=D_i\cap D_j$, we have
\[
\mu(\Pi^*(X),p)=\mu(\wt\fa_X,p).\,\mu(\Pi^*(X),D_i).\,\mu(\Pi^*(X),D_j)\le
\]
\[
\le\mu(\wt\fa_Y,p).\,\mu(\Pi^*(Y),D_i).\,\mu(\Pi^*(Y),D_j)=\mu(\Pi^*(Y),p)
\]
whereas in case II, we have
\[
\mu(\Pi^*(X),p)=\mu(\wt\fa_X,p).\,\mu(\Pi^*(X),D_i)\le\mu(\wt\fa_Y,p).\,\mu(\Pi^*(Y),D_i)=\mu(\Pi^*(Y),p)\,.\qed
\]

Let us finish the proof of theorem \ref{t:2}. Let $p$ be a non-rational singularity of $\wt\fa_X$. By claim \ref{cl:32}, after a formal change of variables $\Phi$ and multiplication by a constant, we have
\[
\Phi^*(\Pi^*(X))=x^p\,y^q\,\left(x\frac{\pa}{\pa x}+\la\,y\frac{\pa}{\pa y}\right):=\wh{X}\,,
\]
where $p+q\ge1$ and $\la\notin\Q$. Set also $\wh{Y}:=\Phi^*(\Pi^*(Y))$. Note that $[X,Y]=0$ implies that $[\wh{X},\wh{Y}]=0$.

On the other hand, the vector field $Z:=q\,x\frac{\pa}{\pa x}-p\,y\frac{\pa}{\pa y}$ commutes with $\wh{X}$: $[Z,\wh{X}]=0$. Since $\wh{X}\wedge Z\not\equiv0$ we can write
\[
\wh{Y}=\phi\,\wh{X}+\psi\,Z\,\,,\,\phi,\psi\in\wh\M_2\,\implies\,
\]
\[
0=[\wh{X},\wh{Y}]=\wh{X}(\phi).\,\wh{X}+\wh{X}(\psi).\,Z\,\,\implies
\] 
$\phi$ and $\psi$ are (formal) meromorphic first integrals of $\wh{X}$ and also of $x\frac{\pa}{\pa x}+\la\,y\frac{\pa}{\pa y}$. But, since $\la\notin\Q$, $\phi$ and $\psi$ must be constants, which implies $\phi=c\in\C$ and $\psi=0$ (because $\mu(\wh{Y},p)>1$). Hence, $Y=c.\,X$ and $Y\in\C.\,X$.
\qed
\vskip.1in

This motivates the following problem:

\begin{prob}\label{pr:4}
Let $X\in\X_2$ be a non-dicritical G.C with an isolated singularity at $0\in\C^2$ and $DX(0)=0$.
Assume also that $r(X)=2$. Is it true that $X$ has a non-constant holomorphic first integral?
\end{prob}

When $DX(0)\ne0$ then the answer is negative, as shows the following example: 
\begin{ex}
{\rm Let
\[
X=x\frac{\pa}{\pa x}-y\frac{\pa}{\pa y}+x\,y\,\left(x\frac{\pa}{\pa x}+\la\,y\frac{\pa}{\pa y}\right)
\]
where $\la\notin\Q$. Then $C(X)=\left<X,x\frac{\pa}{\pa x}-y\frac{\pa}{\pa y}\right>_\C$, so that $r(X)=d(X)=2$, but $\I(X)=\C$.}
\end{ex}

\subsection{Homogeneous and quasi-homogeneous vector fields.}\label{ss:34}
Let $P=(p_1,...,p_n)\in\Z_{\ge1}^n$ with $gcd(p_1,...,p_n)=1$. We say that $f\in\O_n$ is $P$
quasi-homogenous of degree $k\in\N$ if
\begin{equation}\label{eq:1}
f(t^{p_1}.\,z_1,...,t^{p_n}.\,z_n)=t^k.\,f(z_1,...,z_n)\,,\,\forall\,z=(z_1,...,z_n)\in(\C^n,0)\,.
\end{equation}
It is known that (\ref{eq:1}) implies that $f$ is a polynomial.
We would like to observe the following: let $S$ be the linear vector field $S=\sum_{j=1}^np_j\,z_j\,\frac{\pa}{\pa z_j}$. Then
\[
\text{(\ref{eq:1})}\,\iff\,S(f)=k.\,f\,.
\]
We will say also that $f$ is $S$ quasi-homogeneous.
When $S=R=\sum_jz_j\,\frac{\pa}{\pa z_j}$, the radial vector field, then $R(f)=k.\,f$ $\iff$ $f$ is homogeneous of degree $k$.

By analogy, in the case of vector fields, we will say that $X$ is $S$ quasi-homogeneous if $[S,X]=k.\,X$ for some $k\in\Z$. For instance, if $S=R$ then $[R,X]=k.\,X$ $\iff$ the coefficients of $X$ are homogeneous polynomials of degree $k+1$. 

Note also thar, if $S=\sum_jp_j\,z_j\frac{\pa}{\pa z_j}$ and $X=\sum_jX_j\frac{\pa}{\pa z_j}$ then $[S,X]=k\,X$ $\iff$ $X_j$ is $P$ quasi-homogeneous of degree
$k+p_j$, $1\le j\le n$.

\begin{rem}\label{r:21}
{\rm Let $S=\sum_{j=1}^np_j\,z_j\frac{\pa}{\pa z_j}$, $p_j\in\N$, $1\le j\le n$. 
Let $F_k:=\{f\in\O_n\,|\,S(f)=k.\,f\}$ and $E_k:\{X\in\X_n\,|\,[S,X]=k.\,X\}$. 
The following facts are well known:
\begin{itemize}
\item[1.] If $f\in F_k$ (resp. $X\in E_k$) then $f$ is a polynomial (resp. $X$ is a polynomial vector field).
\item[2.] $\O_n=\bigoplus_{k\ge0}F_k$ and $\X_n=\bigoplus_{k\in\Z}E_k$. In particular, any $f\in\O_n$
(resp. any $X\in\X_n$) can be expressed as a convergent series $f=\sum_{k\ge0}f_k$ (resp.
$X=\sum_{k\in\Z}X_k$), where $f_k\in F_k$ (resp. $X_k\in E_k$), for all $k$.
Analogously, any $\wh{f}\in\wh\O_n$ (resp. $\wh{X}\in\wh\X_n$) can be decomposed as a formal power series in
$\bigoplus_{k\ge0}F_k$ (resp. $\bigoplus_{k\in\Z}E_k$).
\item[3.] If $f\in F_k$, $g\in F_\ell$ and $X\in E_m$ then $f.\,g\in F_{k+\ell}$, $f.\,X\in E_{k+m}$ and $X(f)\in F_{k+m}$.
\item[4.] Let $f\in E_k\setminus\{0\}$, $k>0$, so that $f(0)=0$. Let $f=\Pi_{j=1}^rf_j^{m_j}$be the decomposition of $f$ into irreducible factors. Then we can assume that $f_j$ is $S$ quasi-homogeneous; $S(f_j)=k_j.\,f_j$, $1\le j\le r$, where $k=\sum_jm_j\,k_j$.
\item[5.] If $X\in E_k$ and $Y\in E_\ell$ then $[X,Y]\in E_{k+\ell}$
\end{itemize}}
\end{rem}

In the two dimensional case we have the following: 

\begin{prop}\label{p:5}
Let $S=p\,x\frac{\pa}{\pa x}+q\,y\frac{\pa}{\pa y}$ and $X=X_1\frac{\pa}{\pa x}+X_2\frac{\pa}{\pa y}\in \X_2$ with $[S,X]=k\,X$. Let
$\om:=i_Xdx\wedge dy=X_1dy-X_2dx$. Assume that $S\wedge X=h.\,\frac{\pa}{\pa x}\wedge\frac{\pa}{\pa y}$, where
$h=p\,x.\,X_2-q.\,y\,X_1\not\equiv0$. Then $h\in F_{k+tr(S)}$, $tr(S)=p+q$, and is an integrating factor of $\om$: $d\left(\frac{\om}{h}\right)=0$.
Moreover, if $h=\Pi_{j=1}^rh_j^{k_j}$ is the decomposition of $h$ into $S$ quasi-homogeneous factors, then
\begin{itemize}
\item[(a).] There exist $\la_1,...,\la_r\in\C$ and $\var$, $S$ quasi-homogeneous, such that
\begin{equation}\label{eq:oh}
\frac{\om}{h}=\sum_{j=1}^r\la_j\,\frac{dh_j}{h_j}+d\left(\frac{\var}{h_1^{k_1-1}...h_r^{k_r-1}}\right)\,.
\end{equation}
\item[(b).] The curves $\Ga_j:=(h_j=0)$, $1\le j\le r$, are all separatrices of $X$ and $S$.
\item[(c).] If $X$ has other separatrices than the $\Ga_{j`s}$ then $X$ is dicritical.
\end{itemize}
\end{prop}

{\it Proof.}
Set $\mu=\frac{\pa}{\pa x}\wedge\frac{\pa}{\pa y}$. If $Y=Y_1\frac{\pa}{\pa x}+Y_2\frac{\pa}{\pa y}\in\X_2$ then $L_Y\mu=-\nabla Y.\,\mu$ where $\nabla Y=\frac{\pa Y_1}{\pa x}+\frac{\pa Y_2}{\pa y}$.
In particular
\[
k\,h\,\mu=k\,S\wedge X=[S,S]\wedge X+S\wedge[S,X]=L_SS\wedge X=L_S(h\,\mu)=
\]
\[
=S(h)\,\mu+h\,L_S\mu=(S(h)-tr(S)\,h)\mu\,\implies\,S(h)=(k+tr(S))\,h\,\implies\,h\in F_{k+tr(S)}\,.
\]
The proof that $h$ is an integrating factor of $\om$ and of (\ref{eq:oh}) can be found in \cite{cema}.

Let us prove that $(h=0)$ is $X$-invariant. 
It is enough to prove that $h|X(h)$ (see \cite{cema}). We have
\[
0=L_XS\wedge X=L_X(h\,\mu)=X(h).\,\mu-h.\,\nabla X.\,\mu\,\implies\,X(h)=\nabla X.\,h\,\implies\,h|X(h)\,.
\]
From the above relation and $h=\Pi_jh_j^{k_j}$ we get
\[
\nabla X=\frac{X(h)}{h}=\sum_{j=1}^rk_j\,\frac{X(h_j)}{h_j}\,\implies\,h_j|X(h_j)\,,\,1\le j\le r\,\implies
\]
$\Ga_j$ is $X$-invariant, $1\le j\le r$.

Let us prove (c). Let $t\in\C\mapsto S_t:=exp(t.\,S)$ be the flow of $S$. We assert that $S_t$ sends separatrices of $X$ onto separatrices of $X$.

In fact, relation $[S,X]=k\,X$ is equivalent to $S_t^*(X)=e^{k\,t}.\,X$.
Let $f$ be an equation of a separatrix of $X$, so that $X(f)=g.\,f$. If $t\in\C$ is fixed, we have
\[
S_t^*(X(f))=S^*_t(X)(S_t^*(f))=e^{k\,t}.\,X(f\circ S_t)=
\]
\[
=S_t^*(g.\,f)=g\circ S_t.\,f\circ S_t\,\implies\,f\circ S_t|X(f\circ S_t)\,\implies
\]
$(f\circ S_t=0)$ is a separatrix of $X$.

If $X$ has a separatrix $\Ga\notin\{\Ga_1,...,\Ga_r\}$ then this separatrix cannot be
$S_t$-invariant, because otherwise its equation would be contained in the equation $h=0$
(remember that $S\wedge X=h\,\mu$). In this case, the set $\{S_t(\Ga)\,|\,t\in\C\}$ would contain a non-countable set of separatrices of $X$ and $X$ is dicritical.
\qed
\vskip.1in
Concerning the existence of non constant first integrals in the case of dimension two, we have the following:

\begin{prop}\label{p:6}
Let $X,Y\in\X_2\setminus\{0\}$ be $S$ quasi-homogeneous, where $S=p\,x\frac{\pa}{\pa x}+q\,y\frac{\pa}{\pa y}$, $p,q\in\Z_{\ge1}$. Assume that $Y\in C(X)$, $S\wedge X=h\,\mu$ and $X\wedge Y=f\,\mu$, where $h\not\equiv0$.
Then $f/h$ is a first integral of $X$. Moreover, if $X\in E_k$, $Y\in E_\ell$ and $0\ne \ell\ne k$ then $f/h$ is non-constant.
\end{prop}

{\it Proof.}
We assume $f\ne0$.
As we have seen $S\wedge X=h\,\mu$ implies that $X(h)=\nabla X.\, h$. On the other hand, $[X,Y]=0$ implies that
\[
0=L_XX\wedge Y=L_X(f.\,\mu)=(X(f)-\nabla X.\,f)\mu\,\implies\,X(f)=\nabla X.\,f\,\implies
\]
\[
\frac{X(f)}{f}-\frac{X(h)}{h}=0\,\implies\,X\left(\frac{f}{h}\right)=0\,\implies
\]
$f/h$ is a first integral of $X$.

Assume now that $X\in E_k$, $Y\in E_\ell$ and $0\ne\ell\ne k$. We have two possibilities:

{1$^{st}$.} $f\ne0$. We have
\[
L_SX\wedge Y=(k+\ell)\,X\wedge Y=L_S(f\,\mu)=S(f)\,\mu-f\,tr(S)\,\mu\,\implies\,f\in F_{k+\ell+tr(S)}\,.
\]
Similarly, $h\in F_{k+tr(S)}$. In particular, $S(f/h)=\ell\,f/h\ne0$ and so $f/h$ is a non constant first integral of $X$.

{2$^{nd}$.} $f=0$. In this case, there exists a vector field $Z$, with isolated singularity at $0$, and $\phi,\psi\in\O_2\setminus\{0\}$, such that $X=\phi.\,Z$ and $Y=\psi.\,Z$, so that $Y=g.\,X$, $g=\psi/\phi$. In particular, $[X,Y]=0$ implies that $X(g)=0$. Since $X\in E_k$ and $Y\in E_\ell$,
$\ell\ne k$, the function $g$ cannot be constant.
\qed

\begin{thm}\label{t:3}
Let $X\in\X_2$ be $S=p\,x\frac{\pa}{\pa x}+q\,y\frac{\pa}{\pa y}$ quasi-homogeneous with $[S,X]=k.\,X$, $k>0$. Assume that 
$X$ has an isolated singularity at $0\in\C^2$ and is non-dicritical.
If $C(X)\ne\C.\,X$ then $\I(X)=\C\{f\}$, where $f\ne0$ is $S$ quasi-homogeneous and $C(X)=\C\{f\}.\,X$.
\end{thm}

\begin{rem}
{\rm The condition $[S,X]=k.\,X$, $k>0$, implies that $DX(0)$ is nilpotent.}
\end{rem}

{\it Proof.}
The hypothesis implies that there exists $Y\in C(X)$ such that $[X,Y]=0$ and $Y\notin\C.\,X$.
Our purpose is to prove that $Y=\var(f).\,X$, where $f\in\I(X)$ and $\var\in\O_1$.

Let $Y=\sum_{j\in\Z}Y_j$ be the decomposition of $Y$ into $\X_2=\bigoplus_{j\in\Z}E_j$. Note that
\[
\sum_{j\in\Z}[X,Y_j]=[X,Y]=0\,\implies\,[X,Y_\ell]=0\,,\,\forall\ell\in\Z\,,
\]
because $[X,Y_\ell]\in E_{k+\ell}$, $\forall\ell$.
Let us see how $Y_\ell$ looks like. 
\vskip.1in
$1^{st}$ case: $\ell\notin\{0,k\}$. Let $Y_\ell\ne0$ be such that $[X,Y_\ell]=0$. As we have seen in proposition \ref{p:4}, in this case $X$ has a non-constant meromorphic first integral, say $f$. Since $X$ is non-dicritical we can assume that $f$ is holomorphic and minimal, so that $\I(X)=\C\{f\}$: if $f$ was pure meromorphic then $X$ would have infinitely many separatrices and would be dicritical.
We can assume that $f$ is quasi-homogeneous.

In fact, if $f=\sum_{j\ge1}f_j$ is the decomposition of $f$ into $\O_2=\bigoplus_{j\ge0}F_j$ then
$X(f)=0$ implies that $X(f_j)=0$ for all $j\ge1$. Since $f\ne0$, there is $r$ such that $f_r\ne0$, so that $f=f_r\in F_r$ because $f$ is minimal.
We are going to prove that $Y_\ell(f)=0$.

Assume that $Y_\ell(f)\ne0$. Relation $[X,Y_\ell]=0$ implies that
\[
X(Y_\ell(f))=Y_\ell(X(f))=0\,\implies\,Y_\ell(f)\in\C\{f\}\,\implies
\]
\[
Y_\ell(f)=\sum_{i\ge1}a_i\,f^i\,,\,a_i\in\C\,.
\]
Since $Y_\ell(f)\in F_{\ell+r}$ and $a_j.\,f^j\in F_{jr}$ we must have
\begin{equation}\label{eq:6}
Y_\ell(f)=a_m\,f^m\,\,\text{and $a_j=0$, $\forall j\ne m$.}
\end{equation}
In this case, we have also $r\,m=\ell+r$ and $\ell=r(m-1)$ $\implies$ $m>1$.

Relation (\ref{eq:3}) implies also that
\begin{equation}\label{eq:4}
Y_\ell(f^{1-m})=(1-m)\,a_m\,.
\end{equation}
On the other hand, $S(f)=r\,f$ implies that $f^{m-1}\,S(f)=r\,f^m$, and so
\[
\left(Y_\ell-\frac{a_m}{r}\,f^{m-1}\,S\right)(f)=0\,\implies\,\left(Y_\ell-\frac{a_m}{r}\,f^{m-1}\,S\right)\wedge X=0\,.
\]
Since $0$ is an isolated singularity of $X$, we get
\begin{equation}\label{eq:5}
Y_\ell=\frac{a_m.\,f^{m-1}}{r}\,S+g.\,X\,,\,g\in\O_2\,.
\end{equation}
Note that $g\in F_{\ell-k}$, because $Y_\ell,\,g.\,X\in E_\ell$ and $X\in E_k$. 

Relation (\ref{eq:5}) implies also:
\[
\left.
\begin{matrix}
X\wedge Y_\ell=\la_1\,f^{m-1}\,S\wedge X&,\,\la_1=-a_m/r\\
S\wedge Y_\ell=g\,S\wedge X&\,\\
\end{matrix}
\right\}
\implies
\,\frac{f^{m-1}}{g}\in\M(Y_\ell)\,,
\]
where in the last implication we have used proposition \ref{p:4}.
From this and (\ref{eq:4}) we get
\[
Y_\ell\left(\frac{1}{g}\right)=Y_\ell\left(\frac{f^{m-1}}{g}.\,f^{1-m}\right)=(1-m)\,a_m\,.
\]
However, since $g\in F_{\ell-k}$ and $Y_\ell\in E_\ell$ we get
\[
S\left(Y_\ell\left(\frac{1}{g}\right)\right)=k.\,Y_\ell\left(\frac{1}{g}\right)\,\implies\,a_m=0\,\implies\,Y_\ell(f)=0\,,
\]
as we wished. Finally, since $f\in\I(X)\cap\I(Y_\ell)$ we get $Y_\ell=h.\,X$, so that $X(h)=0$ and $h\in\C\{f\}$. In particular, we have $Y_\ell\in\I(X).\,X$.
\vskip.1in
$2^{nd}$ case: $\ell=k$. In this case, we have shown in proposition \ref{p:4} that if $X\wedge Y_k=\psi\,\mu$ and $S\wedge X=\phi\,\mu$, where $\mu=\frac{\pa}{\pa x}\wedge\frac{\pa}{\pa y}$, then $X(\psi/\phi)=0$ (note that $\phi\ne0$, because otherwise $X$ would be a multiple of $S$). Since $L_S\mu=-tr(S)\,\mu$ we get $\psi\in F_{2k+tr(s)}$ and $\phi\in F_{k+tr(s)}$ and so $S(\psi/\phi)=k.\,\psi/\phi$.
We have two possibilities:

{2.1.} $\psi/\phi=c$, where $c\in\C$. In this case, necessarily $c=0$ and $\psi=0$, which implies $Y_k=g\,X$, where $g\in F_0$, so that $g$ is a constant and $Y_k\in\C.\,X$.

{2.2.} $\psi/\phi$ is not a constant. In this case, $X$ has a non constant first integral. Since $X$ is non-dicritical, it has a holomorphic minimal first integral, say $f$, where $f\in F_k$. As in the
$1^{st}$ case, $Y_k(f)\in\C\{f\}$ and $Y_k(f)=a_m\,f^m$, where $f\in F_r$ and $k=r(m-1)$, so that $m>1$.
As in the $1^{st}$ case, we can write
\[
Y_k=\frac{a_m}{r}\,f^{m-1}\,S+g\,X\,,
\]
but now $g\in F_0$ and so it is a constant. Hence,
\[
0=[X,Y_k]=\frac{a_m}{r}\,f^{m-1}\,[X,S]=-\frac{a_m\,k}{r}\,f^{m-1}\,S\,\implies\,a_m=0\,\implies\,Y_k\in\C.X\,.
\]
Therefore, in both cases we get $Y_k\in\C.\,X\sub\I(X).\,X$.
\vskip.1in
{$3^{rd}$ case: $\ell=0$.} We assert that $Y_0=0$. In fact, as before, set $X\wedge Y_0=\psi\,\mu$ and $S\wedge X=\phi\,\mu$,
$\phi\ne0$, so that $X(\psi/\phi)=0$. Note that $\psi,\phi\in F_{k+tr(s)}$ and so $S(\psi/\phi)=0$. Therefore, if $\psi/\phi$ is not a constant then $X$ and $S$ would be colinear, a contradiction: $X$ would be dicritical.
Hence, $X\wedge Y_0=c\,S\wedge X$ $\implies$ $X\wedge(Y_0+c\,S)=0$ $\implies$ $Y_0=-c\,S$, because $X$ has an isolated singularity at $0$. If
$c\ne0$ then $X\wedge S=0$, a contradiction. Therefore $Y_0=0$.  

This finishes the proof of theorem \ref{t:2}.
\qed

\subsection{The dicritical case.}\label{ss:dic}
In this section we will assume that $X\in\X_2$ is dicritical. In this case, if
$\Pi\colon(M,D)\to(\C^2,0)$ is Seidenberg's reduction of the singularities of $X$, then some of the irreducible components of the divisor $D$ are dicritical: non-invariant for the strict transform
$\wt\fa_X$ of the foliation $\fa_X$, defined by $X$.

For instance, the foliation whoose leaves are the levels of a non-constant holomorphic function is non-dicritical, whereas the foliation whoose leaves are the levels of of pure meromorphic function is dicritical.

\begin{ex}\label{ex:3}
{\rm An interesting dicritical foliation (see \cite{cema}) is the one whoose leaves are the levels of the meromorphic function $(y^2+x^3)/x^2$. The associated vector field is
\[
X=2x\,y\frac{\pa}{\pa x}-(2y^2-x^3)\frac{\pa}{\pa y}\,.
\]
As a consequence of the next result we will prove that $C(X)=\C.\,X$.}
\end{ex}

\begin{ex}\label{ex:4}
{\rm An example of dicritical vector field $X\in\X_2$ for which $\mu(X,0)=n\ge2$ and $r(X)=2$ is
\[
X=x^n\frac{\pa}{\pa x}+y^n\frac{\pa}{\pa y}
\]
If $n\ge2$ then $C(X)=\C.\,x^2\frac{\pa}{\pa x}+\C.\,y^2\frac{\pa}{\pa y}$ and $r(X)=d(X)=2$.
If $n=1$ then $X$ is the radial vector field and $d(X)=4$.

When $n\ge2$ then $X$ has the meromorphic first integral $(y^{n-1}-x^{n-1})/x^{n-1}.\,y^{n-1}$.}
\end{ex}

In the next result we consider a germ of dicritical vector field $X\in\X_2$ with an isolated singularity at $0\in\C^2$.
Let
\begin{equation}\label{eq:re}
(M_r,E_r)\overset{\Pi_r}\longrightarrow...(M_2,E_2)\overset{\Pi_2}\longrightarrow(M_1,E_1)\overset{\Pi_1}\longrightarrow(M_0,E_0)=(\C^2,0)
\end{equation}
be the blowing-up process of resolution of singularities of $X$.
Denote by $\wt\fa_m$ the strict transform of $\fa_X$ by the composition
$\wt\Pi_m:=\Pi_1\circ...\circ\Pi_m\colon(M_m,E_m)\to(\C^2,0)$.

\begin{definition}
{\rm Let $Z\in\X_2$ with first non-zero jet $Z_\nu=j^{\nu}(Z,0)$, $\nu\ge1$, such that $Z_\nu\wedge R\equiv0$. We will say that $Z$ has a purely radial singularity at $0$ if $\nu=1$. In this case we have $Z=\a.\,R+h.o.t.$, where $R$ is the radial vector field and $\a\in\C^*$.

If $0\in\C^2$ is an isolated singularity of $Z$ and $\nu>1$ we will say that $Z$ has non purely radial singularity at $0$ (briefly n.p.r.s).
In this case, we have necessarily
\[
Z=f.\,R+h.o.t.\,,
\]
where $f$ is a homogeneous polynomial of degree $\nu-1$.

}
\end{definition}

We will assume that at some step of the resolution, say $\Pi_s\colon(M_s,E_s)\to(M_{s-1},E_{s-1})$, $1\le s\le r$, we have the following:
\begin{itemize}
\item[(1).] $\Pi_s$ is the blowing-up at a dicritical singularity $p_o\in E_{s-1}$ of $\wt\fa_{s-1}$.
\item[(2).] If $s>1$ then $p_o$ belongs
\begin{itemize}
\item[(2.a).] either to a unique divisor $\wt{D}\sub E_{s-1}$, which is a non-dicritical divisor of $\wt\fa_{s-1}$,
\item[(2.b).] or to a corner $\wt{D_1}\cap\wt{D_2}$ of $E_{s-1}$, where both divisors are non-dicritical for $\wt\fa_{s-1}$.
\end{itemize}
\item[(3).] The germ of $\wt\fa_{s-1}$ is defined by a germ at $p_o$ of vector field $\wt{X}$ with a n.p.r.s at $p_o$. In other words, in some coordinates $(x,y)$ around $p_o$, the first non-zero jet $\wt{X}_\nu$ of $\wt{X}$ is of the form
\begin{equation}\label{eq:fr}
\wt{X}_\nu=f.\,R\,,
\end{equation}
where $R=x\frac{\pa}{\pa x}+y\frac{\pa}{\pa y}$ is the radial vector field and $f=f(x,y)$ is a non-constant homogeneous polynomial of degree $\nu-1\ge1$.
\end{itemize}

\begin{thm}\label{t:4}
If $X$ has an isolated singularity at $0\in\C^2$ and is as above then $C(X)=\C.\,X$.
\end{thm}

{\it Proof.}
The proof will be based in the following lemma:
\begin{lemma}\label{l:33}
Let $Z,W\in\X_2$ be such that $[W,Z]=0$ and $\Pi\colon(M,E)\to(\C^2,0)$ be a blowing-up process, where $E$ is the exceptional divisor.
Let $\wt\fa_Z$ and $\wt\fa_W$ denote the strict transforms of the foliations $\fa_Z$ and $\fa_W$ by $\Pi$, respectively.
Suppose that the exceptional divisor $E$ has an irreducible component $D$ such that:
\begin{itemize}
\item[(a).] $D$ is dicritical for $\wt\fa_Z$.
\item[(b).] $\Pi^*(W)|_D\equiv0$.   
\end{itemize}
Then $Z\wedge W\equiv0$.
\end{lemma}

{\it Proof.}
Since $D$ is dicritical for $\wt\fa_Z$, if we fix a generic point $p\in D$ then there are local coordinates $(U,(x,y))$ around $p$ such that
\begin{itemize}
\item[(i).] $x(p)=y(p)=0$ and $E\cap U=D\cap U=(y=0)$.
\item[(ii).] $\wt\fa_Z$ is transverse to $D$ at $p$.
\end{itemize}
In particular, the germ of $\wt\fa_Z$ is represented by a vector field $\wt{Z}$ transverse to $D$ at $p$.
After a local change of variables and taking a smaller $U$ if necessary, we can assume that $\wt{Z}|_U=\frac{\pa}{\pa y}$.
Since $\Pi|_{U\setminus D}$ is a biholomorphism, we must have $\Pi^*(Z)|_U=\var\frac{\pa}{\pa y}$, where $\var(q)\ne0$ for all
$q\in U\setminus D$.

Now, let
\[
\Pi^*(W)|_U=A(x,y)\frac{\pa}{\pa x}+B(x,y)\frac{\pa}{\pa y}\,,
\]
where $A,B\in\O(U)$. From $[Z,W]=0$ we get
\[
\left[\var\frac{\pa}{\pa y},A(x,y)\frac{\pa}{\pa x}+B(x,y)\frac{\pa}{\pa y}\right]=[\Pi^*(Z),\Pi^*(W)]=0
\]
A direct computation shows that the component of $\frac{\pa}{\pa x}$ of $[\Pi^*(Z),\Pi^*(W)]$ is $\pm\,\var.\,A_y$, which implies that $A_y=0$ and $A=A(x)$, so that $\Pi^*(W)=A(x)\frac{\pa}{\pa x}+B(x,y)\frac{\pa}{\pa y}$. Since $\Pi^*(W)|_D\equiv0$ we get $A\equiv0$, and so $\Pi^*(Z)\wedge\Pi^*(W)=0$ on $U$. But this implies $Z\wedge W\equiv0$.
\qed
\vskip.1in
As an application we will prove theorem \ref{t:4} when $p_o=0\in\C^2$, that is when $X=\sum_{j\ge\nu}X_j$, where $X_\nu=f.\,R$, $f$ homogeneous of degree $\nu-1\ge1$.
\begin{cor}\label{c:33}
If $X$ has a non purely radial singularity at $0\in\C^2$ and $0$ is an isolated singularity of $X$ then $C(X)=\C.\,X$.
\end{cor}

{\it Proof.}
Fix $Y\in C(X)$.
Let $\Pi\colon(M,D)\to(\C^2,0)$ be the blow-up at $0\in\C^2$.
Then $D$ is dicritical for $\Pi^*(X)$. We will divide the proof in two cases:

{\it $1^{st}$ case:} $\mu(Y,0)\ge2$. In this case $\Pi^*(Y)|_D\equiv0$ and we can apply lemma \ref{l:33} to show that
$Y\wedge X\equiv0$. Since $X$ has an isolated singularity at $0\in\C^2$ we get $Y=h.\,X$, where $h\in\I(X)$. But since $X$ is dicritical we get $h\in\C$ and $Y\in\C.\,X$.

{\it $2^{nd}$ case:} $\mu(Y,0)=1$. We will see that this is impossible.
In fact, in this case we have $X\wedge Y=h.\,\frac{\pa}{\pa x}\wedge\frac{\pa}{\pa y}$, where $h(0)=0$ and $h\not\equiv0$.
We have seen that any irreducible component of $h$ is invariant by both vector fields $X$ and $Y$.
Let $g$ be an irreducible component of $h$. By lemma \ref{l:31} there exists $\la\in\C$ such that $Y-\la X=g.\,Z$, where $Z\in\X_2$.
Since $\mu(Y,0)=1$ and $\mu(X,0)\ge2$ we get
\[
\mu(g.\,Z)=\mu(Y-\la.\,X,0)=1\,,\,\implies\,\mu(g,0)+\mu(Z,0)=1\,\implies
\]
$\mu(g,0)=1$ and $\mu(Z,0)=0$, because $g(0)=0$ and $g\not\equiv0$. Since $Z(0)\ne0$, after a change of variables we can suppose that $Z=\frac{\pa}{\pa x}$, so that $Y-\la.\,X=g.\,\frac{\pa}{\pa x}$, where $\mu(g,0)=1$. Let $g_1=a\,x+b\,y$ be the linear part of $g$ at $0\in\C^2$. We assert that $a=0$ and $b\ne0$. 

In fact, if $a\ne0$ then for some $\la'\ne\la$ the origin will be a saddle-node of $W:=Y-\la'\,X$, but this is impossible: corollary \ref{c:31} implies that $W$ is holomorphically normalizable and in \S\,\ref{ss:22} it is proved that the pencil generated by $X$ and $Y$ must be equivalent to
\[
\rho\mapsto Z_\rho=\rho.\,y\frac{\pa}{\pa y}+x^p\left(x\frac{\pa}{\pa x}+\ep\,y\frac{\pa}{\pa y}\right)\,,
\]
but then $\ep=1$ and $X=Z_0=x^p\,R$ and $0$ is not an isolated singularity of $X$.
Hence, $a=0$ and after dision by $b$ we can assume that $Y-\la\,X=y\frac{\pa}{\pa x}$.
However, again by \S\,\ref{ss:22}, we must have $X\in \C\{y\}.\left<R,\frac{\pa}{\pa x}\right>$ which implies that $X$ cannot have an isolated singularity at $0\in\C^2$ with $X_\nu=f.\,R$, $deg(f)\ge1$.
\qed
\vskip.1in
In the general case, the idea is similar. 
Recall the blowing-up process of the resolution of singularities of $X$ in (\ref{eq:re}) with $r\ge1$ steps. 
In the $k^{th}$ step $\wt\Pi_k\colon(M_k,E_k)\to(\C^2,0)$ we have called $\wt\fa_k$ the foliation induced by the strict transform $\wt\Pi_k^*(X)$. From now on we will assume $r\ge2$ and that the point $p_o$ which is a n.p.r.s of $\wt\fa_{s-1}$ appears in the $(s-1)$-step, where $s\ge2$.
\vskip.1in
Fix $Y\in C(X)$. Given $1\le m\le r$ let us denote by $\wt\fa_m^Y$ the strict transform of the foliation defined by $Y^*:=\wt\Pi_m^*(Y)$.
In order to apply lemma \ref{l:33}, we have to prove that it is possible to find a n.p.r.s $p_o\in Sing(\wt\fa_{s-1})$ with the property that $\mu(Y^*,p_o)\ge2$. In this case, if we blow-up at $p_o$, $\Pi_s\colon(M_s,D)\to(M_{s-1},p_o)$ then $D$ will be dicritical for $\wt\fa_s$ and $\Pi_s(Y^*)|_D\equiv0$, so that we can apply lemma \ref{l:33}.
\vskip.1in
In order to simplify the proof we will assume the following about the blowing-up process:
\begin{itemize}
\item{} When we pass from the $(m-1)^{th}$ step to the $m^{th}$ step by $\Pi_m\colon(M_m,E_m)\to(M_{m-1},E_{m-1})$ we don't blow-up at a point
$q\in Sing(\wt\fa_{m-1})$ if it is a n.p.r.s or if it is purely radial.
In other words, the blow-up $\Pi_m$ is done at a point $p\in Sing(\wt\fa_{m-1})$ only if it is not a simple singularity and if $\Pi^{-1}(p)=D\sub E_m$ is a non-dicritical divisor for $\wt\fa_m$.
\end{itemize}
Although the final foliation in this process has non-simple singularities, if $\wt\Pi_n=\Pi_1\circ\Pi_2\circ...\circ\Pi_n\colon(M_n,E_n)\to(\C^2,0)$
is the final step, with this convention, then $\wt\fa_n$ satisfies the following:
\begin{itemize}
\item[1.] All irreducible components of $E_n$ are non-dicritical for $\wt\fa_n$.
\item[2.] A non-simple singularity $p\in Sing(\wt\fa_n)$ is, either purely radial, or n.p.r.s.
\item[3.] $\wt\fa_n$ has at least one n.p.r.s, say $p_o\in E_n$.
\end{itemize}

If we suppose that $0\in\C^2$ is not n.p.r.s then:

\begin{claim}\label{cl:34}
We have two possibilities for $X$:
\begin{itemize}
\item[$1^{st}$.] $\mu(X,0)\ge2$. In this case, in any step $\wt\Pi_m\colon(M_m,E_m)\to(\C^2,0)$ of the blowing-up process, in all non-dicritical divisors $D\sub E_m$ we have $\wt\Pi^*(X)|_D\equiv0$.
\item[$2^{nd}$.] $\mu(X,0)=1$ and the linear part of $X$ at $0$ is nilpotent. Moreover, if $D$ is a non-dicritical irreducible component of $E_m$ which is not the strict transform of $E_1$ by
\[
\Pi_2\circ...\circ\Pi_m\colon(M_m,E_m)\to(M_1,E_1)
\]
then $\wt\Pi_m^*(X)|_D\equiv0$.
\end{itemize}
\end{claim}

{\it Proof.}
The $1^{st}$ assertion is proved applying inductively lemma \ref{l:32} in the process.
We leave the details for the reader.

On the other hand, if $\mu(X,0)=1$, the first blow-up is not dicritical for $X$ and $\la_1$ and $\la_2$ are the eigenvalues of $DX(0)$ then:
\begin{itemize}
\item[(i).] If $\la_1,\la_2\ne0$ then in the process of resolution of $X$ there is no n.p.r.s.
\item[(ii).] If $\la_1\ne0$ and $\la_2=0$ then $0$ is a saddle-node for $X$, which is not possible with our hypothesis. 
\end{itemize}
Therefore, $DX(0)$ is nilpotent and we can assume that $DX(0)=y\frac{\pa}{\pa x}$.
In this case, in the first blow-up $\Pi_1\colon(M_1,E_1)\to(\C^2,0)$ then $\Pi_1^*(X)$ has only one singularity at $E_1\simeq\p^1$ corresponding to the direction $p=(y=0)$. This singularity is of algebraic multiplicity $\mu(\Pi_1^*(X),p)=2$.
Therefore, the $2^{nd}$ assertion is also consequence of lemma \ref{l:32}.
\qed
\vskip.1in
Now, given $Y\in C(X)$ let $\wt\fa_n^Y$ be as before, the strict transform of the foliation associated to $\wt\Pi_n^*(Y)$.
Suppose, by contradiction that $X\wedge Y\not\equiv0$.

Note first that all irreducible components of $E_n$ are non-dicritical for $\wt\fa_n^Y$.
Let us prove this fact.

Suppose first that $\mu(X,0)\ge2$. In this case, if $D$ is an irreducible component of $E_n$ then $\wt\Pi_n^*(X)|_D\equiv0$ by claim \ref{cl:34}.
On the other hand, if $D$ was dicritical for $\wt\fa_n^Y$ then $X\wedge Y\equiv0$ by lemma \ref{l:33}, a contradiction.

Suppose now that $\mu(X,0)=1$ and $DX(0)=y\frac{\pa}{\pa x}$. Since $\wt\Pi_n^*(X)|_D\equiv0$ except for $D_1$, the strict transform of $E_1$, we conclude from lemma \ref{l:33} that the unique irreducible component that could be dicritical for $\wt\fa_n^Y$ is $D_1$.
On the other hand, if $D_1$ was dicritical for $\wt\fa_n^Y$ then necessarily $E_1$ is dicritical for the first blow-up of $Y$.
This implies that the first non-zero jet of $Y$ is of the form $Y_k=g.\,R$, where $g$ is a homogeneous polynomial. The polynomial $g$ is necessarily non constant, for otherwise by Poincar\'e linearization theorem we can assume that $Y=R$, which implies that $X$ is a linear vector field contradicting the hypothesis. Now, with an argument similar to the $2^{nd}$ case in the argument of the proof of corollary \ref{c:33} it can be proved that this is impossible.
We leave the details for the reader.
\qed
\vskip.1in
Given an irreducible component $D\sub E_m$ we denote as $\mu(Y,D)$ (resp. $\mu(X,D)$) the order of annihilation of $\wt\Pi_m^*(Y)$ (resp. $\wt\Pi_m^*(X)$) along $D$. In other words, given $p\in D$ and a local coordinate system $(U,(x,y))$ such that $D\cap U=(y=0)$ then $\mu(Y,D)=k$ if $\wt\Pi_m^*(Y)=y^k\,\wt{Y}$ where $y\nmid\,\wt{Y}$.
We have the following:
\begin{claim}\label{cl:36}
Let $p\in D\sub E_m$. Assume that:
\begin{itemize}
\item[(a).] $p$ is not a singularity of $\wt\fa_m^Y$.
\item[(b).] $\mu(Y,D)=k\ge0$.
\end{itemize}
Then:
\begin{itemize}
\item[(1).] If $k=0$ then $p$ is not a singularity of $\wt\fa_m$.
\item[(2).] If $k\ge1$ and $p\in Sing(\wt\fa_m)$ then $p$ is a non-degenerate singularity of $\wt\fa_m$. Moreover, the germ of $\wt\fa_m$ at $p$ is equivalent to the Poincar\'e-Dulac normal form $(k\,x+\a.\,y^k)\frac{\pa}{\pa x}+y\frac{\pa}{\pa y}$.
\end{itemize}
\end{claim}

{\it Proof.}
Since $D$ is non-dicritical for $\wt\fa_m^Y$ and $p\notin Sing(\wt\fa_m^Y)$ we can find a local coordinate system $(U,(x,y))$ around $p$ such that $D\cap U=(y=0)$ and  $\wt\fa_m^Y$ is defined by $\wt{Y}=\frac{\pa}{\pa x}$.
Since $\mu(Y,D)=k$, we can assume that $\wt\Pi_m^*(Y)|_U=y^k.\,\frac{\pa}{\pa x}$.
Let $\Pi^*(X)|_U=y^\ell\,\wt{X}$, where $\ell\ge0$, $\wt{X}=A(x,y)\frac{\pa}{\pa x}+B(x,y)\frac{\pa}{\pa y}$ 
and $y\nmid\wt{X}$. From $[\wt\Pi_m^*(X)\,,\,\wt\Pi_m^*(Y)]=0$ we get
\[
0=\left[y^k.\,\frac{\pa}{\pa x}\,,\,y^\ell\,\wt{X}\right]=y^\ell.\,\left[y^k.\,\frac{\pa}{\pa x}\,,\,\wt{X}\right]\,\implies
\]
\[
\left[y^k.\,\frac{\pa}{\pa x}\,,\,A(x,y)\frac{\pa}{\pa x}+B(x,y)\frac{\pa}{\pa y}\right]=0\,.
\]
As the reader can check directly, the last relation implies that $B_x=0$ and $y.\,A_x=k.\,B$, so that

\[
\left\{
\begin{matrix}
\wt{X}=A(y)\frac{\pa}{\pa x}+B(y)\frac{\pa}{\pa y}&\text{, if $k=0$.}\\
\wt{X}=\left(a(y)+k\,b(y)\,x\right)\frac{\pa}{\pa x}+y\,b(y)\frac{\pa}{\pa y}&\text{, if $k\ge1$.}\\
\end{matrix}
\right.
\]

If $k=0$ then, either $A(0)\ne0$, or $B(0)\ne0$, because $y\nmid\wt{X}$ and so $\wt{X}(p)\ne0$.

If $k\ge1$ and $p\in Sing(\wt\fa_m)$ then $a(0)=0$ and $b(0)\ne0$, for otherwise $y|\wt{X}$.
In particular, the eigenvalues of $D\wt{X}(0)$ are $b(0)$ and $k\,b(0)$.
This implies the last assertion of (2).
\qed
\vskip.1in
\begin{cor}\label{c:34}
If in some step of the resolution of $X$, say $\Pi_m\colon(M_m,E_m)\to(M_{m-1},E_{m-1})$, where $1\le m<n$, we explode at a point $p$ that is not a singularity of $\wt\fa_{m-1}^Y$ then $\Pi_{m}\circ...\circ\Pi_n(p_o)\ne p$, where $p_o\in E_n$ is the n.p.r.s singularity of $\wt\fa_n$.
\end{cor}

{\it Proof.}
The proof is by contradiction: if not, then let $m$ be the smallest step in which we explode at a point $p\in E_{m-1}$ which is not a singularity of $\wt\fa_{m-1}^Y$. Then, since $p$ is a singularity of $\wt\fa_{m-1}$, by claim \ref{cl:36} the germ of $\wt\fa_{m-1}$ at $p$ is represented by a vector field $\wt{X}$ which is equivalent to the Poincar\'e-Dulac normal form
\[
\wt{X}=(k\,x+\a.\,y^k)\frac{\pa}{\pa x}+y\frac{\pa}{\pa y}\,.
\]
But when we continue the process after the resolution of $\wt{X}$ we don't obtain any n.p.r.s singularity, a contradiction.
\qed
\vskip.1in

Let us see how looks like $\wt\Pi_m^*(Y)$ in a neighborhood of a point $p\in E_m$, $1\le m\le n$.
Denote as $\wt\fa_m^Y$ the strict transform of $\fa_Y$ by $\wt\Pi_m$.
For the first blow-up we have the following:

\begin{claim}\label{cl:35}
\begin{itemize}
\item[(a).] If $\mu(X,0)\ge2$ and $\mu(Y,0)=1$ then $DY(0)$ is nilpotent.
\item[(b).] If $\mu(X,0)=1$ and $DX(0)$ is nilpotent then $DY(0)$ is also nilpotent and $DY(0)\wedge DX(0)=0$.
\end{itemize} 
\end{claim}

{\it Proof.}
Let $DY(0)=S+N$, where $S=\la_1\,x\frac{\pa}{\pa x}+\la_2\,y\frac{\pa}{\pa y}$ is semi-simple and $N$ is nilpotent and $[S,N]=0$.
We assert that $S=0$.

In fact, suppose by contradiction that $S\ne0$. Let $\mu(X,0)=k$ and $X_k=j^k(X,0)$ be the first non-zero jet of $X$.
Note that
\[
[S+N,X_k]=0\,\,\implies\,\,[S,X_k]=[N,X_k]=0\,\,.
\]
Suppose first that $k\ge2$. In this case, we must have $\la_1,\la_2\ne0$, for otherwise for $\a\in\C^*$ we have $Z:=Y+\a.\,X\in C(X)$ has a saddle-node at $0\in\C^2$, which is not possible with our hypothesis.
On the other hand, if $\la_1,\la_2\ne0$, then $[S,X_k]=0$ implies that $S$ has a resonance and necessarily, after multiplying by a constant is equivalent to either $S=x\frac{\pa}{\pa x}+n\,y\frac{\pa}{\pa y}$, $n\in\N$, or to $S=m\,x\frac{\pa}{\pa x}-n\,y\frac{\pa}{\pa y}$, $m,n\in\N$.

Let $p_o\in E_n$ be the n.p.r.s singularity of $\wt\fa_n$ and consider the sequence of images of $p_o$
\[
p_o\overset{\Pi_n}\longrightarrow p_1\overset{\Pi_{n-1}}\longrightarrow p_2\overset{\Pi_{n-2}}\longrightarrow...\overset{\Pi_1}\longrightarrow p_n=0\in\C^2
\]
By corollary \ref{c:34} all points in the sequence are singularities of the strict transform of $\fa_Y$.
Given $1\le j\le n$ let $D_j:=\Pi^{-1}_j(p_{n-j+1})$ be the irreducible component of $E_j$ obtained in the blowing-up $\Pi_j$.
We assert that $\mu(\wt\Pi_j^*(Y),p_{n-j})=1$ and $\mu(\,\wt\Pi^*_j(Y),D_j)=0$, $1\le j\le n$.

The above assertion is consequence of the following: 

A. When we blow-up a non-degenerated and non radial singularity of a germ at $0\in\C^2$ of vector field $Z$, say by $\Pi\colon(\wt\C^2,D)\to(\C^2,0)$, $\mu(\Pi^*(Z),D)=0$ and we have two possibilites for the singularities of $\Pi^*(Z)$:
\begin{itemize}
\item[A.1.] $\Pi^*(Z)$ has two non-degenerated singularities in the divisor $D$.
\item[A.2.] $\Pi^*(Z)$ has only one singularity in $D$, which is saddle-node. This happens only when $DZ(0)$ is equivalent to $R+y\frac{\pa}{\pa x}$.
\end{itemize}

B. When we blow-up at a saddle-node then it appears two singularities at the divisor, one non-degenerated and the other a saddle-node.

If we apply A and B inductivelly we obtain the assertions. Moreover, at the end of the process the vector field $Y^*:=\wt\Pi_n^*(Y)$ has, either a non-degenerated singularity, or a saddle-node at $p_o$. However this is not possible because $\wt\Pi_n^*(X)$ has a n.p.r.s singularity at $p_o$, as the reader can check.
\vskip.1in
Suppose now that $k=1$. In this case $DX(0)$ is nilpotent and we can assume that $DX(0)=y\frac{\pa}{\pa x}$.
Since $DY(0)$ and $DX(0)$ commute we can write $DY(0)=a\,R+b\,DX(0)$, $R$ the radial vector field.
We assert that $a=0$.
In fact, if $a\ne0$ then by Poincar\'e's linearization theorem we can assume that $Y=a.\,R+b\,DX(0)$.
However, in this case $X$ cannot satisfy hypothesis (3) (see \S\,\ref{ss:22}). 
Therefore, $DY(0)$ is nilpotent and $DX(0)\wedge DY(0)=0$.
\qed

\vskip.1in
Now, let $p_o\in E_n$ be a n.p.r.s: the germ at $p_o$ of $\wt\fa_n$ is defined by $\wt{X}$, where the first non-zero jet of $\wt{X}$ is $\wt{X}_k=f.\,R$, $f$ homogeneous of degree $\nu-1\ge1$.
In order to apply lemma \ref{l:33} it is suficient to prove that $\mu(\wt\Pi_n^*(Y),p_o)\ge2$.
Note that, in principle, $p_o$ could be a pole of $\wt\Pi_n^*(Y)$, if at some steps of the process we had blow-up at points that are not singularities of strict transform of $\fa_Y$.
However, this doesn't happens by corollary \ref{c:34}.
\vskip.1in
{\it Proof that $\mu(\wt\Pi_n^*(Y),p_o)\ge2$.} Suppose by contradition that $\mu(\wt\Pi_n^*(Y),p_o)=1$.
Recall that $p_o\in Sing(\wt\fa_n^Y)$ by corollary \ref{c:34}.
This implies that $p_o$ is an isolated singularity of $\wt\Pi_n^*(Y):=\wt{Y}$, so that $\wt{Y}$ defines the germ at $p_o$ of $\wt\fa_n^Y$.

Choose coordinates $(U,(x,y))$ such that $x(p_o)=y(p_o)=0$ and
\begin{itemize}
\item[(iii).] $p_o\in D$, where $D$ is an irreducible component of $E_n$. Moreover, $D\cap U=(y=0)$ and $\mu(\wt\Pi_n^*(X),D)=k\ge0$.
\item[(iv).] The germ of $\wt\fa_n$ is defined by $\wt{X}\in\X_2$ where $j^\nu(\wt{X})=X_\nu=f.\,R$, where $f$ is homogeneous of degree $\nu-1$.
Moreover, $\wt\Pi_n^*(X)=y^k\,\wt{X}$.
\item[(v).] $\wt{Y}=A(x,y)\frac{\pa}{\pa x}+B(x,y)\frac{\pa}{\pa y}$, where $D\wt{Y}(0)\ne0$. We have seen before that $D\wt{Y}(0)$ is nilpotent, so that we can suppose that $D\wt{Y}(0)=y\frac{\pa}{\pa x}$.
\end{itemize}
Since $[\wt\Pi_n^*(X),\wt{Y}]=0$, if we set $\wt\Pi_n^*(X)\wedge\wt{Y}=g.\,\frac{\pa}{\pa x}\wedge\frac{\pa}{\pa y}$, then any irreducible component of $g$ is $\wt{Y}$ and $\wt\Pi_n^*(X)$ invariant.
Let $h$ be an irreducible component of $g$. Then there exists $\la\in\C$ such that $h|\wt{Y}-\la.\,\wt\Pi_n^*(X)$:
\[
\wt{Y}-\la.\,\wt\Pi_n^*(X)=h\,Z\,.
\]
From the above relation we get
\[
\mu(h\,Z,0)=\mu(h,0)+\mu(Z,0)=1\,\,\implies\,\,\mu(h,0)=1\,\,\text{and}\,\,\mu(Z,0)=0\,\,.  
\]
As in the proof of corollary \ref{c:33} we can assume that $Z=\frac{\pa}{\pa x}$ (recall that $D$ is non-dicritical for $\wt{Y}$) and $h=y$.
The vector field $\wt\Pi_n^*(X)$ in the new coordinate system can still be written as $y^k\wt{X}$ where $0=p_o$ is an isolated singularity of $\wt{X}$, as the reader can check.
In particular,
\[
0=\left[y\frac{\pa}{\pa x}\,,\,y^k\wt{X}\right]=y^k\,\left[y\frac{\pa}{\pa x}\,,\,\wt{X}\right]\,\implies\,\left[y\frac{\pa}{\pa x}\,,\,\wt{X}\right]=0\,\implies
\]
\[
\wt{X}\in\,\C\{y\}.\left<R,\frac{\pa}{\pa x}\right>\,\,,
\]
but then $p_o$ cannot be a singularity of the type n.p.r.s.
This finishes the proof of theorem \ref{t:4}.
\qed
\vskip.1in
We have the following consequence of theorem \ref{t:4}:
\begin{cor}
Let $X\in\X_2$ be a dicritical vector field with an isolated singularity at $0\in\C^2$ such that $r(X)=2$.
Let $\Pi\colon(M,E)\to(\C^2,0)$ be the minimal resolution of the singularities of $X$ and $D_i\sub E$ be a dicritical divisor of $E$.
Then $D_i^2=-1$ (the self intersection number) and $D_i$ is obtained by the blowing-up of a purely radial singularity on a non-dicritical divisor of the previous step of the resolution.
\end{cor}

Another interesting result is the following:

\begin{prop}\label{p:7}
Let $X\in\X_2$ be a dicritical vector field with an isolated singularity at $0\in\C^2$.
Asume that $DX(0)$ is nilpotent. If $r(X)=2$ then $d(X)=2$.
\end{prop}

{\it Proof.}
As in the proof of theorem \ref{t:4}, consider the blowing-up process $\wt\Pi_r=\Pi_1\circ...\circ\Pi_r\colon(M_r,E_r)\to(\C^2,0)$ as in \ref{eq:re}.
By theorem \ref{t:4} during the process it never appears a n.p.r.s singularity.
As a consequence, there is at least a dicritical exceptional divisor $D\sub E_r$ such that $D^2=-1$, which was obtained by blowing-up at a purely radial singularity on a non-dicritical divisor of a previous step.
Without lost of generality, we will assume that $D=\Pi_r^{-1}(p_o)$, where $p_o\in D'\sub E_n$, $n=r-1$, is the purely radial singularity, and $D'$ is the irreducible component of $E_{n}$ that contains $p_o$.
Given $1\le m\le n$, set as before $\wt\Pi_m=\Pi_1\circ...\circ\Pi_m\colon(M_m,E_m)\to(\C^2,0)$ and let $\wt\fa_m$ be the strict transform of $\fa_X$ at the level $m$: $\wt\Pi_m=\wt\Pi_m^*(\fa_X)$.

\begin{rem}\label{r:37}
{\rm If $D'\sub E_n$ is the irreducible component that contains $p_o$ then $\mu(\wt\Pi_n^*(X),D')\ge1$.
This is a direct consequence of claim \ref{cl:34}.}
\end{rem}

\begin{claim}\label{cl:37}
There are coordinates $(U,(x,y))$ around $p_o$ such that:
\begin{itemize}
\item[1.] $D'\cap U=(y=0)$ and $x(p_o)=y(p_o)=0$.
\item[2.] $\wt\Pi_n^*(X)=y^k.\,R$, where $k=\mu(\wt\Pi_n^*(X),D')\ge1$ and $R$ is the radial vector field.
\end{itemize}
\end{claim}

{\it Proof.}
Let $(U',(u,v))$ be local coordinates such that $u(p_o)=v(p_o)=0$ and $D'\cap U'=(v=0)$.
Since $\mu(\wt\Pi_n^*(X),D')=k$ we can write $\wt\Pi_n^*(X)=v^k.\,\wt{X}$, where $v\nmid\wt{X}$.
By assumption $D\wt{X}(0)=\a\,R$, where $\a\ne0$, so that by Poincar\'e's linearization theorem, after a change of variables, we can assume that $\wt{X}=\a\,R$, $D'\cap U'=(v=0)$ and $\wt\Pi_n^*(X)=\phi.\,v^k.\,R$, where $\phi(0)\ne0$.
If we consider a change of variables of the form $h(u,v)=(\rho.\,u,\rho.\,v)=(x,y)$, where $\rho(0)\ne0$, then
\[
h_*(\wt\Pi_n^*(X))=\phi.\,\rho^{-k-1}.\,(R(\rho)+\rho).\,y^k\left(x\frac{\pa}{\pa x}+y\frac{\pa}{\pa y}\right)\,.
\]
Now, we use the fact that the differential equation $R(\rho)+\rho=\phi^{-1}\,\rho^{k+1}$ has a solution $\rho$ such that $\rho^k(0)=\phi(0)$.
We leave the details for the reader.
\qed
\vskip.1in
Now, let $Y\in C(X)$ be such that $X\wedge Y\not\equiv0$. 
We can assume that $0$ is an isolated singularity of $Y$: if not then instead of $Y$ take $Y'=Y+\a.\,X$, where $\a\ne0$.

Denote as $\wt\fa_m^Y$ the strict transform of $\wt\Pi_m^*(\fa_Y)$ on $M_m$, where $\wt\fa_0=\fa_Y$.
As before, consider the sequence of images of $p_o$:
\[
p_0=p_o,\,p_1=\Pi_n(p_o),\,p_2=\Pi_{n-1}(p_1),...,p_n=\Pi_1(p_{n-1})=0\in\C^2\,.
\]
Note that $p_{n-m}\in E_m$ if $1\le m\le n$.

Define $S(Y)=\{m\,|\,0\le m\le n$ and $p_{n-m}\in Sing(\wt\fa_m^Y)\}$.

\begin{claim}\label{cl:38}
There exists $0<\ell<n$ such that $S(Y)=\{k\,|\,0\le k\le\ell\}$. Moroever $\mu(\wt\Pi_n^*(Y),D)=1$.
In other words $p_o\notin Sing(\wt\fa_n^Y)$ and $\wt\Pi_n^*(Y)=y.\,\wt{Y}$, where $\wt{Y}(p_o)\ne0$ and $(y=0)$ is a local equation on $D'$.
\end{claim}

{\it Proof.}
Of course $0=p_n\in Sing(\fa_Y)$. We will prove at the end that $S(Y)\ne\{0,...,n\}$. Assuming that $S(Y)\ne\{0,...,n\}$, there exists $0\le m<n$ such that $0,...,m-1\in S(Y)$, but $m\notin S(Y)$.

Let $\wt{Y}$ and $\wt{X}$ be germs at $p_{n-m}$ of vector fields representing $\wt\fa_m^Y$ and $\wt\fa_m$, respectively.
By claim \ref{cl:36}, there are local coordinates $(U,(u,v))$ around $p_{n-m}$ and $\ell\ge1$ such that $u(p_{n-m})=v(p_{n-m})=0$ and
\begin{itemize}
\item[(i).] $E_m\cap U=(v=0)$.
\item[(ii).] $\wt{Y}=\frac{\pa}{\pa u}$ and $\wt\Pi_m^*(Y)=v^\ell.\,\wt{Y}$.
\item[(iii).] $\wt{X}=(\ell.\,u+\a\,v^\ell)\frac{\pa}{\pa u}+v\,\frac{\pa}{\pa v}$.
\end{itemize}
Note that $\a=0$ because otherwise the singularity $p_{n-m}$ would not be dicritical.
Following the resolution of $\ell.\,u\,\frac{\pa}{\pa u}+v\,\frac{\pa}{\pa v}$, we see that the radial singularity $p_o$ appears after $\ell-1$ blowing-ups and $m=n-\ell+1$.
Moreover, the composition $\Pi:=\Pi_{n-\ell+2}\circ...\circ\Pi_n\colon(M_n,p_o)\to(M_{n-\ell+1},p_{n-m})$ in the chart that appears the radial singularity $p_o$ is of the form $\Pi(x,y)=(x.\,y^{\ell-1},y)=(u,v)$ with inverse $\Pi^{-1}(u,v)=(u/v^{\ell-1},v)=(x,y)$.
It follows that
\[
\wt\Pi_n^*(Y)=\Pi^*\wt\Pi_m^*(Y)=\Pi^*\left(v^\ell\frac{\pa}{\pa u}\right)=y\frac{\pa}{\pa x}\,.
\]
This proves (b).
It remains to prove that $S(Y)\ne\{0,...,n\}$.

Suppose by contradiction that $S(Y)\ne\{0,...,n\}$.
In this case, all points in the sequence $p_0,...,p_n$ are singularities of the strict transform of $\fa_Y$.
Note that $DY(0)$ is nilpotent. The proof of this fact is similar to the proof of claim \ref{cl:35}. We leave the details for the reader.

Since all points of the sequence are singularities of the strict transform of $\fa_Y$, when we apply lemma \ref{l:32} inductively we find that $\mu(\wt\Pi_m^*(Y),p_{n-m})\ge2$ if $1\le m\le n$. In particular, $\mu(\wt\Pi_n^*(Y),p_o)\ge2$. This implies that when we apply the blowing-up $\Pi_{n+1}\colon (M_r,D)\to(M_n,p_o)$ then we have $\Pi_{n+1}^*\wt\Pi_n^*(Y)|_D\equiv0$, which is impossible by lemma \ref{l:33}.
\qed
\vskip.1in
Let us finish the proof of proposition \ref{p:7}. 
Let $(U,(x,y))$ be a coordinate system around $p_o$ as in claim \ref{cl:37}:
\begin{itemize}
\item[(I).] $x(p_o)=y(p_o)=0$ and $D'\cap U=(y=0)$.
\item[(II).] $\wt\Pi_n^*(X)=y^k.\,R$.
\end{itemize}

We assert that
\[
C(y^k.\,R)=\left\{(\a.\,x+\be\,y)\frac{\pa}{\pa x}\,+\,g.\,R\,\,|\,\a,\be\in\C\,\,\text{and}\,g\,\text{is homogeneous of degree}\,\,k\right\}\,.
\]

In fact, let $Z\in C(y^k\,R)$ and write $Z=\sum_{j\ge0}Z_j$ where $Z_j$ is homogeneous of degre $j\ge0$.
From $[Z,y^k.\,R]=0$ we get
\[
0=k\,y^{k-1}\,Z(y)\,R+y^k\,[Z,R]=y^{k-1}\left(k\,Z(y)\,R-y\,\sum_{j\ge0}(j-1)\,Z_j\right)\,\implies
\]
\begin{equation}\label{eq:17}
k\,Z(y)\,R=y.\,\sum_{j\ge0}(j-1)\,Z_j\,.
\end{equation}
From (\ref{eq:17}) we get $Z_0=0$ and $Z_j\wedge R=0$ if $j\ge2$.
In particular, we can write $Z=Z_1+h.\,R$, where $h\in\O_2$ and $h(0)=0$.
Note that the term of degree $2$ of the right side of (\ref{eq:17}) vanishes and this implies that $Z_1(y)=0$ so that $Z_1=(\a\,x+\be\,y)\frac{\pa}{\pa x}$ and $[Z_1,y^k\,R]=0$. It follows that
\[
[h\,R,y^k\,R]=0\,\,\implies\,\,R(h)=k\,h\,\,\implies\,\,h\,\text{is homogeneous of degree $k$.}
\]

From the above we can assume that $\Pi^*_n(Y)=y\,\frac{\pa}{\pa x}+c\,y^k\,R$, where $c\in\C$.
Let $W\in C(X)$ with $W\wedge X\ne0$. Since $\wt\Pi_n^*(W)|_{D'}\equiv0$, by the above argument we must have
\[
\wt\Pi_n^*(W)=\a\,y\frac{\pa}{\pa x}+g\,R\,,
\]
where $g$ is homogeneous of degree $k$. Since $W-\a\,Y\in C(X)$ and $\wt\Pi_n^*(W-\a\,Y)\wedge R\equiv0$ we obtain $(W-\a\,Y)\wedge X\equiv0$.
Therefore $W=\a\,Y+\be\,X$, $\be\in\C$.
This finishes the proof of proposition \ref{p:7}.
\qed
\vskip.1in
\begin{ex}\label{ex:p7}
{\rm An example satisfying proposition \ref{p:7} is $X=y\frac{\pa}{\pa x}+x^n\,R$, $R$ the radial vector field.
Note that $X$ is dicritical because admits the meromorphic first integral $\frac{1}{y^n}-\frac{n}{n+1}\left(\frac{x}{y}\right)^{n+1}$.
In this case we have $C(X)=\left<X,y^n\,R\right>$.}
\end{ex}
\vskip.1in
As a consequence of proposition \ref{p:7}, we have the following:

\begin{cor}\label{c:37}
Let $X\in\X_2$ be a dicritical germ of vector field with an isolated singularity at $0\in\C^2$.
Then $d(X)\le4$. Moreover:
\begin{itemize}
\item[(a).] If $d(X)=3$ then $X$ is equivalent, modulo a multiplicative constant, to $x\frac{\pa}{\pa x}+n\,y\frac{\pa}{\pa y}$ where $n\in\N\setminus\{1\}$.
\item[(b).] If $d(X)=4$ then $X$ is equivalent, modulo a multiplicative constant, to the radial vector field.
\end{itemize}
\end{cor}

{\it Proof.}
By proposition \ref{p:7} if $r(X)=2$ and $d(X)>2$ then $DX(0)$ cannot be nilpotent.
Since $X$ is dicritical and cannot have a saddle-node at the origin $0\in\C^2$, we obtain that the semi-simple part of $DX(0)$, after multiplication by a constant, is of the form $m\,x\,\frac{\pa}{\pa x}+n\,y\,\frac{\pa}{\pa y}$, where $m,n\in\N$ and $gcd(m,n)=1$.
In fact, from Poincar\'e's linearization and Poincar\'e-Dulac theorems we can assume that, after multiplication by a constant, that $X$ is conjugated to one of the following vector fields:
\begin{itemize}
\item[1.] $m\,x\frac{\pa}{\pa x}+n\,y\frac{\pa}{\pa y}$, where $1<m<n$. In this case we have $d(X)=2$ (see \S\,\ref{ss:22}).
\item[2.] $x\frac{\pa}{\pa y}+n\,y\frac{\pa}{\pa y}$, $n>1$. In this case we have $d(X)=3$.
\item[3.] $X=R$, the radial vector field. In this case we have $d(X)=4$.
\end{itemize}
Therefore, the corollary is a direct consequence of the examples in \S\,\ref{ss:22}.
\qed
\vskip.1in
Another consequence of theorem \ref{t:4} and of proposition \ref{p:7} is the following:

\begin{cor}
Let $X\in\X_2$ with an isolated singularity at $0\in\C^2$. If $d(X)>4$ then $X$ has a non-constant holomorphic first integral.
In particular $d(X)=\infty$.
\end{cor}

{\it Proof.}
Note first that $X$ cannot be dicritical by corollary \ref{c:37}, because $d(X)>4$.
Moreover, $d(X)>2\ge r(X)$ and $X$ has a non-constant meromorphic first integral by proposition \ref{p:2}.
In particular, all leaves of $\fa_X$ are closed and since $X$ is non-dicritical it has a non-constant holomorphic first integral by \cite{mm}.
\qed

\begin{rem}
{\rm When the origin is not an isolated singularity of $X\in\X_2$ then $d(X)$ can be arbitrarily large and finite.
An example is $X_k=(x\,y)^k.\,R$, where $R=x\frac{\pa}{\pa x}+y\frac{\pa}{\pa y}$ is the radial vector field.
It can be checked that
\[
C(X_k)=\left<x\frac{\pa}{\pa x}-y\frac{\pa}{\pa y}\,,\,h(x,y).\,R\,|\,h\,\text{is homogeneous of degree $2k$}\right>_\C\,\,,
\]
which has dimension $d(X)=2k+2$ (see also example \ref{ex:5} in \S\,\ref{ss:21}).}

\end{rem}
\vskip.1in
This motivates the following:
\begin{prob}
{\rm Let $X\in\X_2$ with $d(X)=\infty$ and non-isolated singularity at $0\in\C^2$. Is $\I(X)\ne\C$ ?}
\end{prob}
\begin{rem}
{\rm Example \ref{ex:4} with $n\ge2$ has two distinguished separatrices: $(x=0)$ and $(y=0)$. In fact, consider the pencil of commuting vector fields $X_\la=x^n\frac{\pa}{\pa x}+\la\, y^n\frac{\pa}{\pa y}$, $\la\in\C^*$. Then $(x^{n-1}-\la\,y^{n-1})/x^{n-1}\,y^{n-1}$ is a first integral of $X_\la$. In particular, if $\la\ne0,1$ then the unique commom separatrices of $X_1=X$ and $X_\la$ are the curves $(x=0)$ and $(y=0)$.

If $\la\ne0$, when we blow-up once at $0\in\C^2$, $\Pi\colon(\wt\C^2,D)\to(\C^2,0)$, the strict transform $\wt\fa_\la$ of the foliation defined by $\Pi^*(X_\la)$ has n+1 singularities, all non-degenerated, two of them are non-dicritical and the others dicritical. The non-dicritical singularities correspond to the directions $(x=0)$ and $(y=0)$ and don't change with the parameter. The dicritical singularities move along the divisor $D$ with the parameter $\la$. If $n=2$ the dicritical singularity is purely radial, whereas if $n>2$ then at any dicritical singularity the foliation $\wt\fa_\la$ has a meromorphic first integral which in some cordinate system $(u,v)$ is of the form $v^{n-1}/u$, where the local equation of $D$ is $v=0$.
This motivates the following:
\begin{prob}
{\rm Is the above situation general?
More specifically, suppose that $X\in\X_2$ is dicritical, has an isolated singularity at $0\in\C^2$, $\mu(X,0)\ge2$ and $r(X)=2$. Let $Y\in C(X)$ such that $X\wedge Y\ne0$ and consider the pencil $\la\mapsto X_\la:=X+\la\,Y$.

$1^{st}$ question: is there a blowing-up process $\Pi\colon(M,E)\to(\C^2,0)$, where $E$ has an irreducible component $D$ which contains movables non-degenerated singularities with local meromorphic first integrals?

$2^{nd}$ question: does there exists a holomorphic family $\la\mapsto f_\la$ of non-constant meromorphic first integrals: $X_\la(f_\la)=0$?}
\end{prob}}
\end{rem}

\bibliographystyle{amsalpha}

\vskip.2in

D. Cerveau

Universit\'e de Rennes 1

IRMAR, CNRS UMR 6625

F-35042 Rennes Cedex, France

email adress: dominique.cerveau@univ-rennes1.fr

\vskip.2in

A. Lins Neto

IMPA, Est. D. Castorina, 110, 22460-320, Rio de Janeiro, RJ, Brazil

email adress: alcides@impa.br

\end{document}